\documentclass[11pt,twoside]{article}
\usepackage{srcltx}
\usepackage{calrsfs}

{\catcode`\@=11 \gdef\ps@myheadings{\let\@mkboth\@gobbletwo
 \def\@oddhead
          {\vbox{\noindent
                {\small HEAVY TAILS IN MULTI-SERVER QUEUE}
                                         \hfill\rm\thepage\vskip 4pt\hrule}}%
 \def\@oddfoot{}
 \def\@evenhead
          {\vbox{\noindent\rm\thepage\hfill\small
          S. FOSS AND D. KORSHUNOV\vskip 4pt \hrule}}%
 \def\@evenfoot{}\def\chaptermark##1{}\def\sectionmark##1{}%
 \def\subsectionmark##1{}}

\gdef\@begintheorem#1#2{\trivlist \item[\hskip \labelsep{\indent\bf #1\ #2.}]\it}
\gdef\@opargbegintheorem#1#2#3{\trivlist
      \item[\hskip \labelsep{{\bf #1}\ #2\ (#3).}]\it}
\gdef\@endtheorem{\endtrivlist}
 }

\newcounter{remark}
\newcommand\remark
       {\refstepcounter{remark}
                  {R\,e\,m\,a\,r\,k~\arabic{remark}.\ }}
\newcommand\proof{P\,r\,o\,o\,f}
\newcounter{definition}

\newlength{\myskip}
\newtheorem{Theorem}{Theorem}
\newtheorem{Lemma}{Lemma}
\newtheorem{Corollary}{Corollary}

\newcommand\mysection[1]{
              \refstepcounter{section}
              \section*{\normalsize\bf\thesection.~#1}
                         }

\begin{document}
\thispagestyle{empty}

\section*{\Large\bf Heavy Tails in Multi-Server Queue\footnote{Supported
           by EPSRC grant No.~R58765/01, INTAS Project No.~00-265
           and RFBR grant No.~02-01-00358}\\[2mm]
           \normalsize\rm SERGUEI FOSS \hfill foss@ma.hw.ac.uk\\
           {\it Department of Actuarial Mathematics and Statistics,
           School of Mathematical and Computer Sciences,
           Heriot-Watt University, Edinburgh EH14 4AS, Scotland}\\
           DMITRY KORSHUNOV \hfill korshunov@math.nsc.ru\\
           \it Sobolev Institute of Mathematics,
           Novosibirsk 630090, Russia
          }

\vspace{5mm}

{\small{\bf Abstract.} In this paper, the asymptotic behaviour
of the distribution tail of the stationary
waiting time $W$ in the $GI/GI/2$ FCFS queue is studied.
Under subexponential-type assumptions on the service time distribution,
bounds and sharp asymptotics are given for the probability ${\bf P}\{W>x\}$.
We also get asymptotics for the distribution tail
of a stationary two-dimensional workload vector and of
a stationary queue length. These asymptotics
depend heavily on the traffic load.\\[5mm]
{\bf Keywords:} FCFS multi-server queue,
stationary waiting time, large deviations,
long tailed distribution, subexponential distribution.
}

\mysection{Introduction}

It is well known (see, for example, [\ref{P}, \ref{Ver},
\ref{APQ}])
that in the stable single server {\it first-come-first-served}
queue $GI/GI/1$ with typical interarrival time $\tau$ and
typical service time $\sigma$ the tail of stationary waiting time $W$
is related to the service time distribution tail
$\overline B(x)={\bf P}\{\sigma>x\}$ via the equivalence
\begin{eqnarray}\label{W.single}
{\bf P}\{W>x\} &\sim& \frac{1}{{\bf E}\tau-{\bf E}\sigma}
\int_x^\infty \overline B(y)\,dy \quad\mbox{ as }x\to\infty,
\end{eqnarray}
provided the {\it subexponentiality} of the
{\it integrated tail distribution} $B_I$ defined by its tail
\begin{eqnarray*}
\overline B_I (x) &\equiv& \min\Bigl(1,\
\int_x^\infty \overline B(y)\,dy \Bigr),
\ \ x>0.
\end{eqnarray*}
As usual we say that a distribution $G$ on ${\bf R}^+$
is subexponential (belongs to the class $\mathcal S$)
if $\overline{G*G}(x)\sim2\overline G(x)$ as $x\to\infty$.
The converse assertion is also true, that is, the equivalence
(\ref{W.single}) implies the subexponentiality of $B_I$,
see [\ref{P}, Theorem 1] for the case of Poisson arrival stream
and [\ref{K}, Theorem 1] for the general case.

In this paper we consider the $GI/GI/s$ FCFS queue
which goes back to Kiefer and Wolfowitz [\ref{KW55}].
We have $s$ identical servers, i.i.d.\ interarrival times
$\{\tau_n\}$ with finite mean $a={\bf E}\tau_1$,
and i.i.d.\ service times $\{\sigma_n\}$
with finite mean $b={\bf E}\sigma_1$.
The sequences $\{\tau_n\}$ and
$\{\sigma_n\}$ are mutually independent.  The system is assumed
to be {\it stable}, i.e., $\rho\equiv b/a\in(0,s)$.
We are interested in the asymptotic tail behaviour of the stationary
waiting time distribution ${\bf P}\{W>x\}$ as $x\to\infty$.

It was realized recently
(see, for example, existence results for moments in [\ref{S}],
[\ref{SS}]; an asymptotic hypothesis in [\ref{W}];
asymptotic results for fluid queues fed by heavy-tailed
on-off flows in [\ref{BMZ}])
that the heaviness of the stationary waiting time tail
depends substantially on the load $\rho$ in the system.
More precisely, it depends on $\rho$ via the value
of $k\in\{0,1,\ldots,s-1\}$ for which $k\le\rho<k+1$.
In particular, Whitt conjectured that
\begin{eqnarray*}
{\bf P}\{W>x\} &\sim&
\gamma\left(\int_{\eta x}^\infty\overline B(y)dy\right)^{s-k}
\quad \mbox{ as } x\to\infty,
\end{eqnarray*}
``where $\gamma$ and $\eta$ are positive constants
(as functions of $x$)'' [sic, [\ref{W}]].
In the present paper we show that, in general,
the tail behaviour of $W$ is more complicated.

Let $R(w)=(R_1(w),\ldots,R_s(w))$ be the operator on
${\bf R}^s$ which orders the coordinates of $w\in{\bf R}^s$
in ascending order, i.e., $R_1(w)\le\cdots\le R_s(w)$.
Then the residual work vector $W_n=(W_{n1},\ldots,W_{ns})$
which the $n$th customer observes just upon its arrival
satisfies the celebrated Kiefer--Wolfowitz recursion:
$W_1=i\cdot0$,
\begin{eqnarray*}
W_{n+1} &=& R((W_{n1}+\sigma_n-\tau_{n+1})^+,(W_{n2}-\tau_{n+1})^+,
\ldots, (W_{ns}-\tau_{n+1})^+)\\
&=& R(W_n+e_1\sigma_n-i\tau_{n+1})^+,
\end{eqnarray*}
here $e_1=(1,0,\ldots,0)$, $i=(1,\ldots,1)$ and
$w^+=(\max(0,w_1),\ldots,\max(0,w_s))$.
The value of $W_{n1}$ is the delay which customer $n$ experiences.
In particular, the stationary waiting time $W$
is a weak limit for $W_{n1}$.

The process $W_n$ is a Markov chain in ${\bf R}^s$.
It is well known that, for general multi-dimensional Markov
chains, large deviation problems are very difficult
to solve even for stationary distributions. Usually they can be solved
in low dimensions only, 2 or 3 at most,
see [\ref{IMSh}, \ref{BMumn}].
Almost all known results are derived for so-called Cram\'er
case which corresponds to light-tailed distributions
of jumps. In the heavy-tailed case almost nothing is known
for general multi-dimensional Markov chains.

The process $W_n$ presents a particular but very important
example of a Markov chain in ${\bf R}^s$,
even if we are interested in the first component $W_{n1}$.
As follows from our analysis, the case $s=2$ can be
treated in detail.
The stability condition for this particular case is $b<2a$.
One of the following cases can occur:

(i) the maximal stability case when $b<a$;

(ii) the intermediate case when $b=a$;

(iii) the minimal stability case when $b\in(a,2a)$.

We find the exact asymptotics for ${\bf P}\{W>x\}$
in the maximal and minimal cases.
We also describe the most probable way
for the occurrence of large deviations.
In the intermediate case, we only provide
upper and lower bounds. Then we study the asymptotics
for the tail of the distribution of a stationary two-dimensional workload
vector and give comments on the tail asymptotics of the
stationary queue length.

For $s>2$, the stability condition is $b<sa$.
We hope that, for $s>2$, direct modifications of our arguments may
lead to exact asymptotics in two particular cases when either
$b<a$ (the maximal stability) or $b \in ((s-1)a, sa)$ (the minimal
stability). However, one has to overcome many extra technicalities
for that. Insofar as the case $b\in[a,(s-1)a]$ is concerned,
we are extremely sceptical on the possibility to get any sharp tail
asymptotics in explicit form.

For the two-server queue,
in the maximal stability case, we prove the following:

\begin{Theorem}\label{th.2.max}
Let $s=2$ and $b<a$. When the integrated tail distribution $B_I$
is subexponential, the tail of the stationary waiting time
satisfies the asymptotic relation, as $x\to\infty$,
\begin{eqnarray*}
{\bf P}\{W>x\} &\sim&
\frac{1}{a(2a-b)}\Bigl[(\overline B_I (x))^2
+b\int_0^\infty \overline B_I(x+ya)\overline B(x+y(a-b))dy\Bigr].
\end{eqnarray*}
\end{Theorem}

The proof follows by combining the lower bound given in Theorem
\ref{th.2.max.lower} (Section \ref{sec.2.max.lower})
and the
upper bound given in Theorem \ref{th.2.max.upper}
(Section \ref{sec.2.max.upper}).
Simpler lower and upper bounds for
${\bf P}\{W>x\}$ are given in the following

\begin{Corollary}\label{cor.2.max}
Under the conditions of Theorem {\rm\ref{th.2.max},}
\begin{eqnarray*}
\frac{2a+b}{2a^2(2a-b)} &\le&
\liminf_{x\to\infty}
\frac{{\bf P}\{W>x\}}{(\overline B_I (x))^2}
\le \limsup_{x\to\infty}
\frac{{\bf P}\{W>x\}}{(\overline B_I (x))^2}
\le \frac{1}{2a(a-b)}.
\end{eqnarray*}
\end{Corollary}

In our opinion, in Theorem \ref{th.2.max} it is
possible to obtain a compact expression
for the tail asymptotics of ${\bf P} (W>x)$ only
in the regularly varying case. A distribution $G$
(or its tail $\overline{G}$) is
{\it regularly varying} at infinity
with index $\gamma >0$ (belongs to the class ${\cal RV}$),
if $\overline{G}(x)>0$ for all $x$ and, for any fixed $c>0$,
$\overline{G}(cx)/\overline{G}(x)\to c^{-\gamma}$
as $x\to\infty$.

\begin{Corollary}\label{cor.2.max.reg}
Let $b<a$ and the tail distribution $\overline B$ of service time
be regularly varying with index $\gamma>1$. Then
\begin{eqnarray*}
{\bf P}\{W>x\} &\sim& c'(\overline B_I(x))^2,
\end{eqnarray*}
where
\begin{eqnarray*}
c' &=& \frac{1}{a(2a-b)}\Bigl[1+\frac{b}{\gamma-1}
\int_0^\infty\frac{dz}{(1+za)^{\gamma-1}(1+z(a-b))^\gamma}\Bigr].
\end{eqnarray*}
\end{Corollary}

Recall definitions of a number of classes of heavy-tailed
distributions.
A distribution $G$ is {\it long-tailed} (belongs to the class
${\mathcal L}$) if $\overline{G}(x)>0$ for all $x$ and, for any fixed $t$,
\begin{eqnarray*}
\frac{\overline G(x+t)}{\overline G(x)}
&\to& 1 \ \ \mbox{ as }x\to\infty.
\end{eqnarray*}
A distribution $G$ belongs to the class ${\cal IRV}$ of {\it
intermediate regularly varying distributions} if
$\overline{G}(x) >0$ for all $x$ and
$$
\lim_{c\downarrow 1}\liminf_{x\to\infty}
\frac{\overline{G}(cx)}{\overline{G}(x)} =1.
$$
Clearly, ${\cal RV} \subset {\cal IRV}$.

In the minimal stability case, we prove the following

\begin{Theorem}\label{th.2.min.reg}
Let $s=2$ and $a<b<2a$, $B\in{\mathcal S}$ and
$B_I\in{\cal IRV}$. Then
\begin{eqnarray*}
{\bf P}\{W>x\}
&\sim& \frac{1}{2a-b} \overline B_I\Bigl(\frac{b}{b-a}x\Bigr)
\quad\mbox{ as }x\to\infty.
\end{eqnarray*}
\end{Theorem}

The proof is given in Section \ref{min.stab.proof} and
is based on the lower and upper bounds stated in Theorems
\ref{th.2.min.lower} and \ref{min_upper} respectively.

One can provide simple  sufficient conditions
for $B\in{\mathcal S}$ and  $B_I\in{\cal IRV}$.
Let ${\cal D}$ be the class of all distributions $G$ on
${\bf R}^+$ such that $\overline{G}(x) >0$ for all $x$ and
$\liminf_{x\to\infty}\overline G(2x)/\overline G(x)>0$.
Then the following are known:
(i) ${\cal RV } \subset {\cal IRV} \subset
({\cal L} \bigcap {\cal D})\subset{\mathcal S}$;
(ii) if $B\in {\cal D}$ has a finite first moment,
then $B_I\in{\cal IRV}$ (see e.g. [\ref{BoF}]).
Therefore, if $B\in {\mathcal L} \bigcap {\cal D}$
and has a finite first moment, then $B$ satisfies
the conditions of Theorem \ref{th.2.min.reg}.
Note that the converse is not true, in general:
there exists a distribution $B\in{\mathcal S}$
with a finite first moment such that $B_I\in{\cal IRV}$,
but $B\notin{\mathcal L} \bigcap{\cal D}$
(see Example 2 in [\ref{DFK}, Section 6]).

The paper is organized as follows.
Section \ref{preliminaries} contains some auxiliary results.
In Section \ref{sec.2.max.lower}, we formulate and prove
a result concerning a lower bound for ${\bf P}\{W>x\}$ in
the maximal stability case. The corresponding upper bound
is given in Section \ref{sec.2.max.upper}.
Sections \ref{sec.2.min.lower}, \ref{sec.2.min.upper},
and \ref{min.stab.proof} deal, respectively,
with lower bounds, upper bounds, and asymptotics for
${\bf P}\{W>x\}$ in the minimal stability case.
In Section \ref{workload}, we prove further results related
to the joint distribution of a stationary workload vector.
Comments on the asymptotics for a stationary queue length
distribution may be found in Section \ref{stat.q.l}.

A number of upper and lower bounds for
${\bf P}\{W>x\}$ in $s$-server queue are proposed
in Remarks \ref{rem.s.min}, \ref{rem.s.min.up},
\ref{sserverlower}, and \ref{sserverupper}.

\mysection{Preliminaries\label{preliminaries}}

{\bf \ref{preliminaries}.1. Reduction to deterministic input
stream case in assertions associated with upper bounds.}
Consider a general $GI/GI/s$ queue.
Take any $a'\in(b/s,a)$.
Consider an auxiliary $D/GI/s$ system with the same
service times $\{\sigma_n\}$ and deterministic interarrival
times $\tau_n'\equiv a'$: $W'_1=0$ and
\begin{eqnarray*}
W'_{n+1} &=& R(W_n'+e_1\sigma_n-ia')^+.
\end{eqnarray*}
Let $W'$ be a stationary waiting time in this auxiliary system.

\begin{Lemma}\label{upper.bound}
If ${\bf P}\{W'>x\} \le \overline G(x)$
for some long-tailed distribution $G$, then
\begin{eqnarray*}
\limsup_{x\to\infty} \frac{{\bf P}\{W>x\}}{\overline G(x)} &\le& 1.
\end{eqnarray*}
\end{Lemma}

\proof. Denote $\xi_n=a'-\tau_n$. Put $M_0 = 0$ and,
for $n\ge1$,
\begin{eqnarray*}
M_n &=& \max\{0,\ \xi_n,\
\xi_n+\xi_{n-1},\ \ldots,\ \xi_n+\cdots+\xi_1\}\\
&=& \max (0, \xi_n+M_{n-1}) = (\xi_n+M_{n-1})^+.
\end{eqnarray*}
First, we use induction to prove the inequality
\begin{eqnarray}\label{upper.via.M}
W_n &\le& W'_n+iM_n\quad\mbox{ a.s.}
\end{eqnarray}
Indeed, for $n=1$ we have $0\le0+iM_1$.
Assume
the inequality is proved for some $n$; we prove it for $n+1$.
Indeed,
\begin{eqnarray*}
W_{n+1} &=& R(W_n+e_1\sigma_n-i\tau_{n+1})^+\\
&\le& R(W_n'+iM_n+e_1\sigma_n-i\tau_{n+1})^+\\
&=& R(W_n'+e_1\sigma_n-ia'+i(M_n+\xi_{n+1}))^+.
\end{eqnarray*}
Since $(u+v)^+\le u^++v^+$,
\begin{eqnarray*}
W_{n+1}
&\le& R(W_n'+e_1\sigma_n-ia')^+ +i(M_n+\xi_{n+1})^+
\equiv W_{n+1}'+iM_{n+1},
\end{eqnarray*}
and the proof of (\ref{upper.via.M}) is complete.

Let $M$ be  the weak limit for $M_n$ which exists due to
${\bf E}\xi_1=a'-a<0$ and Strong Law of Large Numbers.
The following stochastic equality holds:
\begin{eqnarray*}
M &=_{\rm st}& \max\{0,\ \xi_1,\
\xi_1+\xi_2,\ \ldots,\ \xi_1+\cdots+\xi_n,\ \ldots\}.
\end{eqnarray*}
Since the random variable $\xi_1$ is bounded from above (by $a'$),
there exists $\beta>0$ such that ${\bf E}e^{\beta\xi_1}=1$.
Then by Cram\'er estimate (see, for example,
[\ref{Cr}, Section 5]), for any $x$,
\begin{eqnarray}\label{exp.bound.for.M}
{\bf P}\{M>x\} &\le& e^{-\beta x}.
\end{eqnarray}
The inequality (\ref{upper.via.M}) implies that
$W \le_{\rm st} W'+M$, where $W'$ and $M$ are independent.
Let a random variable $\eta$ have distribution $G$ and be
independent of $M$. Since $\eta\ge_{\rm st}W'$,
we have $W \le_{\rm st} \eta+M$. Therefore, for any $h>0$,
\begin{eqnarray*}
{\bf P}\{W>x\}
&\le& \int_0^{x-h} {\bf P}\{M>x-y\}{\bf P}\{\eta\in dy\}
+{\bf P}\{\eta>x-h\}\\
&\le& \int_0^{x-h} e^{-\beta(x-y)} G(dy)+\overline G(x-h),
\end{eqnarray*}
by (\ref{exp.bound.for.M}). Integrating by parts yields
\begin{eqnarray*}
\int_0^{x-h} e^{-\beta(x-y)} G(dy)
&=& -e^{-\beta(x-y)} \overline G(y)\Big|_0^{x-h}
+\beta \int_0^{x-h} \overline G(y) e^{-\beta(x-y)}dy\\
&\le& e^{-\beta x} +\beta \int_h^x \overline G(x-y) e^{-\beta y}dy.
\end{eqnarray*}
The distribution $G$ is long-tailed, thus, for any $\varepsilon>0$
there exists $x(\varepsilon)$ such that
\begin{eqnarray*}
\overline G(x-1) &\le& \overline G(x)e^\varepsilon
\end{eqnarray*}
for any $x\ge x(\varepsilon)$. Hence,
there exists $c(\varepsilon)<\infty$ such that
\begin{eqnarray*}
\overline G(x-y) &\le& c(\varepsilon)\overline G(x)e^{\varepsilon y}
\end{eqnarray*}
for any $x\ge x(\varepsilon)$. Take $\varepsilon=\beta/2$. Then
\begin{eqnarray*}
\int_h^x \overline G(x-y) e^{-\beta y}dy
&\le& c(\varepsilon)\overline G(x)\int_h^x e^{-\beta y/2}dy
\le \frac{c(\varepsilon)}{\beta/2}\overline G(x) e^{-\beta h/2}.
\end{eqnarray*}
Hence,
\begin{eqnarray*}
{\bf P}\{W>x\}
&\le& e^{-\beta x} +2c(\varepsilon)\overline G(x) e^{-\beta h/2}
+\overline G(x-h).
\end{eqnarray*}
Taking into account also that $\overline G(x-h)\sim\overline G(x)$
for any fixed $h>0$, we obtain
\begin{eqnarray*}
\limsup_{x\to\infty} \frac{{\bf P}\{W>x\}}{\overline G(x)}
&\le& 2c(\varepsilon) e^{-\beta h/2}+1.
\end{eqnarray*}
Letting $h\to\infty$ yields the conclusion of the Lemma.

{\bf \ref{preliminaries}.2. Reduction to deterministic input
stream case in assertions associated with lower bounds.}
Take any $a'>a$. As in the previous subsection,
consider an auxiliary $D/GI/s$ system with the same service times
$\{\sigma_n\}$ and deterministic interarrival times $\tau_n'\equiv a'$.
Let $W'$ be a stationary waiting time in this auxiliary system.

\begin{Lemma}\label{lower.bound}
If ${\bf P}\{W'>x\} \ge \overline G(x)$
for some long-tailed distribution $G$, then
\begin{eqnarray*}
\liminf_{x\to\infty} \frac{{\bf P}\{W>x\}}{\overline G(x)} &\ge& 1.
\end{eqnarray*}
\end{Lemma}

\proof. Put $\xi_n=\tau_n-a'$, $M_0 = 0$ and
\begin{eqnarray*}
M_n &=& \max\{0,\ \xi_n,\
\xi_n+\xi_{n-1},\ \ldots,\ \xi_n+\cdots+\xi_1\}
=(M_{n-1}+\xi_n)^+.
\end{eqnarray*}
For any $n\ge1$, the following inequality holds:
\begin{eqnarray}\label{lower.via.M}
W_n &\ge& W'_n-iM_n.
\end{eqnarray}
Indeed, by induction arguments,
\begin{eqnarray*}
W_{n+1} &=& R(W_n+e_1\sigma_n-i\tau_{n+1})^+\\
&\ge& R(W_n'-iM_n+e_1\sigma_n-i\tau_{n+1})^+\\
&=& R(W_n'+e_1\sigma_n-ia'-i(M_n+\xi_{n+1}))^+.
\end{eqnarray*}
Since $(u-v)^+\ge u^+-v^+$,
\begin{eqnarray*}
W_{n+1}
&\ge& R(W_n'+e_1\sigma_n-ia')^+ -i(M_n+\xi_{n+1})^+
\equiv W_{n+1}'-iM_{n+1},
\end{eqnarray*}
and the proof of (\ref{lower.via.M}) is complete.

Let $M$ be the weak limit for $M_n$ which exists due to
${\bf E}\xi_1=a-a'<0$ and the Strong Law of Large Numbers.
The inequality (\ref{lower.via.M}) implies that
$W \ge_{\rm st} W'-M$ where $W'$ and $M$ are independent.
Therefore, for any $h>0$,
\begin{eqnarray*}
{\bf P}\{W>x\}
&\ge& {\bf P}\{W'>x+h\} {\bf P}\{M\le h\}
\ge\overline G(x+h) {\bf P}\{M\le h\}.
\end{eqnarray*}
The distribution $G$ is long-tailed,
thus $\overline G(x+h)\sim\overline G(x)$
for any fixed $h>0$ and
\begin{eqnarray*}
\liminf_{x\to\infty} \frac{{\bf P}\{W>x\}}{\overline G(x)}
&\ge& {\bf P}\{M\le h\}.
\end{eqnarray*}
Letting $h\to\infty$, we obtain the desired estimate from below.

{\bf \ref{preliminaries}.3. Adapted versions
of the Law of Large Numbers.}
It is well known that obtaining lower bounds
for systems under assumptions of heavy tails
usually requires
some variant of the Law of Large Numbers. Here we provide
such a tool for the two-server queue.

\begin{Lemma}\label{SLLN.max}
Let $(\xi_n,\eta_n)$, $n=1$, $2$, \ldots, be independent
identically distributed pairs of random variables.
Let the two-dimensional Markov chain
$V_n=(V_{n1},V_{n2})$, $n=1$, $2$, \ldots, be defined
in the following way: $V_1$ has an arbitrary distribution and
\begin{eqnarray*}
V_{n+1} &=& \left\{
\begin{array}{lll}
V_n+(\xi_n,\eta_n),\ &{\rm if}&\ V_{n1}\le V_{n2},\\
V_n+(\eta_n,\xi_n),\ &{\rm if}&\ V_{n1}>V_{n2}.
\end{array}
\right.
\end{eqnarray*}
If ${\bf E}\eta_1<{\bf E}\xi_1$,
then the following convergence in probability holds:
\begin{eqnarray*}
\frac{V_n}{n} &\to&
\Bigl(\frac{{\bf E}\xi_1+{\bf E}\eta_1}{2},\
\frac{{\bf E}\xi_1+{\bf E}\eta_1}{2}\Bigr)
\ \ \mbox{ as }n\to\infty.
\end{eqnarray*}
\end{Lemma}

\proof. Since $V_{n+1,1}+V_{n+1,2}=V_{n1}+V_{n2}+\xi_n+\eta_n$,
by the Law of Large Numbers
\begin{eqnarray}\label{sum.of.VV}
\frac{V_{n1}+V_{n2}}{n} &\to& {\bf E}\xi_1+{\bf E}\eta_1
\quad\mbox{ as }n\to\infty.
\end{eqnarray}

Define a Markov chain $U_n=V_{n2}-V_{n1}$.
If $U_n\ge0$, then $U_{n+1}-U_n=\eta_n-\xi_n$, while
if $U_n<0$, then $U_{n+1}-U_n=\xi_n-\eta_n=-(\eta_n-\xi_n)$,
so, $U_n$ is the oscillating random walk.
Since ${\bf E}\xi_1>{\bf E}\eta_1$, the mean drift
of the chain $U_n$ is negative on the positive half-line
and is positive on the negative half-line.
Therefore, for any sufficiently large $A$,
the set $[-A,A]$ is positive recurrent for this Markov chain.
In particular,
the distributions of $U_n$ are tight.
Hence, $U_n/n \to 0$ in probability as $n\to\infty$.
Together with (\ref{sum.of.VV}),
it implies the desired assertion of Lemma.

The classical Law of Large Numbers
and Lemma \ref{SLLN.max} imply the following

\begin{Corollary}\label{SLLN.max.2}
Let ${\bf E}\eta_1<{\bf E}\xi_1<0$ and $\varepsilon>0$.
Then
\begin{eqnarray*}
{\bf P}\{V_{n1}>0,\,V_{n2}>0\,|\,V_1=(v_1,v_2)\} &\to& 1
\end{eqnarray*}
as $N\to\infty$ uniformly in $n\ge N$
and in $(v_1,v_2)$ on the set
$$
\Bigl\{v_1, v_2>n(|{\bf E}\xi_1|+\varepsilon),\
v_1+v_2>n(|{\bf E}\xi_1|+|{\bf E}\eta_1|+\varepsilon)\Bigr\}.
$$
\end{Corollary}

\begin{Corollary}\label{SLLN.max.3}
Let ${\bf E}\eta_1<{\bf E}\xi_1<0$ and $\varepsilon>0$.
Then
\begin{eqnarray*}
{\bf P}\{V_{n1}>0,\,V_{n2}>0\,|\,V_1=(v_1,v_2)\} &\to& 0
\end{eqnarray*}
as $N\to\infty$ uniformly in $n\ge N$
and in $(v_1,v_2)$ on the complementary set
$$
\overline{\{v_1>n(|{\bf E}\xi_1|-\varepsilon),\
v_2>n(|{\bf E}\xi_1|-\varepsilon),\
v_1+v_2>n(|{\bf E}\xi_1|+|{\bf E}\eta_1|-\varepsilon)\}}.
$$
\end{Corollary}

\begin{Corollary}\label{SLLN.max.4}
Let ${\bf E}\eta_1<0$, ${\bf E}\xi_1>0$,
${\bf E}\eta_1+{\bf E}\xi_1<0$ and $\varepsilon>0$.
Then
\begin{eqnarray*}
{\bf P}\{V_{n1}>x,\,V_{n2}>x\,|\,V_1=(v_1,v_2)\} &\to& 1
\end{eqnarray*}
as $x$, $N\to\infty$ uniformly in $n\ge N$
and in $(v_1,v_2)$ on the set
$$
\Bigl\{v_1>x-n({\bf E}\xi_1-\varepsilon),\
v_2>2x+n(|{\bf E}\xi_1+{\bf E}\eta_1|+\varepsilon)\Bigr\}.
$$
\end{Corollary}

\mysection{The maximal stability case: a lower bound
\label{sec.2.max.lower}}

\begin{Theorem}\label{th.2.max.lower}
Assume $b\in (0,a)$.
Let the integrated service time distribution $B_I$ be long-tailed.
Then the tail of the stationary waiting time $W$ admits
the following estimate from below: as $x\to\infty$,
\begin{equation} \label{firstlb}
{\bf P}\{W>x\}
\ge \frac{1+o(1)}{a(2a-b)} \Biggl[(\overline B_I (x))^2+
b\int_0^\infty \overline B_I (x+ya)\overline B(x+y(a-b))dy\Biggr].
\end{equation}
\end{Theorem}

\remark\label{rem.rou}
From (\ref{firstlb}), one can get the lower bound
in Corollary \ref{cor.2.max}.
Namely, replace $\overline{B}(x+y(a-b))$ by a smaller term
$\overline{B}(x+ya)$
in the integral in the RHS of (\ref{firstlb}). Then the new
integral is equal to $b(\overline{B}_I(x))^2/2a$, and the lower
bound follows since
$$
\frac{1}{a(2a-b)} \Bigl(1+\frac{b}{2a}\Bigr) =
\frac{2a+b}{2a^2(2a-b)}.
$$

\remark\label{rem.s.min}
By use of Strong Law of Large Numbers,
one can get the following result for
$s$-server queue, $s\ge 2$. If $b<a$, then
there exists a constant $K\equiv K(a,b,s)$ such that
$$
{\bf P}\{W>x\} \ge (K+o(1))(\overline{B}_I(x))^s.
$$

We start with some auxiliary results.
The proof of the theorem is given in subsection \ref{sec.2.max.lower}.4.

{\bf \ref{sec.2.max.lower}.1. An integral equality.}

\begin{Lemma}\label{calcul.1}
Let $f(y)$ be an integrable function. Put
$f_I(y)\equiv\int_y^\infty f(z)dz$. Then,
for any positive $\alpha$ and $\beta$, $\alpha>\beta$,
\begin{eqnarray*}
J &\equiv&
\int_0^\infty \int_0^\infty
f(\alpha y{+}\beta z)f(\beta y{+}\alpha z)dydz\\
&& \hspace{30mm} =\ \frac{(f_I(0))^2}{\alpha^2{-}\beta^2}
- \frac{2\beta}{\alpha^2{-}\beta^2}
\int_0^\infty f_I(\alpha u)f(\beta u)du.
\end{eqnarray*}
\end{Lemma}

\proof. Put $u=\alpha y{+}\beta z$ and
$v=\beta y{+}\alpha z$. Then
\begin{eqnarray*}
J &=& \frac{1}{\alpha^2{-}\beta^2}
\int_0^\infty f(u)du
\int_{\beta u/\alpha}^{\alpha u/\beta} f(v)dv\\
&=& \frac{1}{\alpha^2{-}\beta^2}
\int_0^\infty f(u)f_I(\beta u/\alpha)du
-\frac{1}{\alpha^2{-}\beta^2}
\int_0^\infty f(u)f_I(\alpha u/\beta)du\\
&=& \frac{\alpha}{\alpha^2{-}\beta^2}
\int_0^\infty f(\alpha u)f_I(\beta u)du
-\frac{\beta}{\alpha^2{-}\beta^2}
\int_0^\infty f(\beta u)f_I(\alpha u)du.
\end{eqnarray*}
Integration by parts yields
\begin{eqnarray*}
\int_0^\infty f_I(\beta u)f(\alpha u)du
&=& \frac{(f_I(0))^2}{\alpha}-
\frac{\beta}{\alpha}
\int_0^\infty f_I(\alpha u)f(\beta u)du.
\end{eqnarray*}
By substituting this equality into the previous one,
we arrive at the conclusion of the Lemma.

{\bf \ref{sec.2.max.lower}.2. Some calculations
with two big service times.}
Fix $\varepsilon>0$ and put $b'=b-\varepsilon$.
For $k$ and $l$, $k<l\le n$, define the events
$A_{nkl}$ and $C_{nkl}$ by the equalities
\begin{eqnarray*}
A_{nkl} &=&
\Bigl\{\sigma_k>x+(l-k)a+(n-l)(a-b'),\ \sigma_l>x+(n-l)(a-b'),\\
&& \hspace{40mm}
\sigma_k+\sigma_l>2x+(l-k)a+(n-l)(2a-b')\Bigr\}
\end{eqnarray*}
and
\begin{eqnarray*}
C_{nkl} &=& \bigcap_{\stackrel{j=1}{j\ne k,l}}^n\Bigl\{
\sigma_j\le x+(n-j)(a-b')\Bigr\}.
\end{eqnarray*}
Note that the events $A_{nkl}\cap C_{nkl}$
are disjoint for different pairs $(k,l)$.
Due to the existence of ${\bf E}\sigma$,
uniformly in $n\ge1$ and $k<l\le n$,
\begin{eqnarray}\label{slln.for.C}
{\bf P}\{\overline C_{nkl}\}
&\le& \sum_{j=0}^\infty {\bf P}\{\sigma_1>x+j(a-b')\} \to 0
\quad\mbox{ as }x\to\infty.
\end{eqnarray}

\begin{Lemma}\label{int.2}
Assume $b\in (0,a)$.
Let the integrated tail distribution $B_I$ be long-tailed.
Then, for any fixed $N\ge1$ and for any $\varepsilon >0$, as $x\to\infty$,
\begin{eqnarray*}
\lim_{n\to\infty}
\sum_{\stackrel{k,l=1}{k<l}}^{n-N} {\bf P}\{A_{nkl}\}
&\sim& \frac{1}{a(2a-b')} \Biggl[(\overline B_I (x))^2
+ b'\int_0^\infty \overline B_I(x+ya)\overline B(x+y(a-b'))dy
\Biggr].
\end{eqnarray*}
\end{Lemma}

\proof. Put
\begin{eqnarray*}
A_{kl}' &=&
\{\sigma_1>x+ka+l(a{-}b'),\ \sigma_2>x+l(a{-}b'),\
\sigma_1+\sigma_2>2x+ka+l(2a{-}b')\},
\end{eqnarray*}
so that ${\bf P}\{A_{nkl}\}={\bf P}\{A'_{l-k,n-l}\}$ and
\begin{eqnarray}\label{throw.A.prime}
\lim_{n\to\infty}\sum_{\stackrel{k,l=1}{k<l}}^{n-N} {\bf P}\{A_{nkl}\}
&=& \lim_{n\to\infty}\sum_{l=N}^{n-1}\sum_{k=1}^{n-l-1} {\bf P}\{A_{kl}'\}
= \sum_{l=N}^\infty\sum_{k=1}^\infty {\bf P}\{A_{kl}'\}.
\end{eqnarray}
Consider also the events
\begin{eqnarray*}
A(y,z) &=& \{\sigma_1>x{+}ya{+}z(a{-}b'),\
\sigma_2>x{+}z(a{-}b'),\
\sigma_1{+}\sigma_2>2x{+}ya{+}z(2a{-}b')\},
\end{eqnarray*}
which satisfy $A(k,l)=A'_{kl}$.
Since the probability ${\bf P}\{A(y,z)\}$
is non-increasing in $y$ and $z$, we have the inequalities
\begin{eqnarray}\label{diff.i-.i+}
I_-\equiv\int_N^\infty\int_1^\infty {\bf P}\{A(y,z)\}dydz
&\le& \sum_{l=N}^\infty\sum_{k=1}^\infty {\bf P}\{A_{kl}'\}
\nonumber\\
&\le& \int_0^\infty\int_0^\infty {\bf P}\{A(y,z)\}dydz
\equiv I_+.
\end{eqnarray}
The values of integrals $I_-$ and $I_+$ are close
to each other in the following sense:
\begin{eqnarray*}
\lefteqn{I_+-I_-}\\
&\le& \int_0^N\int_0^\infty {\bf P}\{A(y,z)\}dydz
+ \int_0^\infty\int_0^1 {\bf P}\{A(y,z)\}dydz\\
&\le& N{\bf P}\{\sigma_2>x\}
\int_0^\infty{\bf P}\{\sigma_1>x+ya\}dy
+{\bf P}\{\sigma_1>x\}
\int_0^\infty{\bf P}\{\sigma_1>x+z(a{-}b')\}dz.
\end{eqnarray*}
Recall that the distribution  $\overline B_I(x)$ is long tailed,
which is equivalent to $\overline B(x)=o(\overline B_I(x))$.
Therefore, as $x\to\infty$,
\begin{eqnarray*}
I_+-I_- &\le&
\frac{N+1}{a-b'}\overline B(x) \overline B_I(x)
=o((\overline B_I(x))^2).
\end{eqnarray*}
Now it follows from (\ref{diff.i-.i+}) that, as $x\to\infty$,
\begin{eqnarray}\label{integral.repr}
\sum_{l=N}^\infty\sum_{k=1}^\infty {\bf P}\{A_{kl}'\}
&=& \int_0^\infty\int_0^\infty {\bf P}\{A(y,z)\}dydz
+o((\overline B_I(x))^2).
\end{eqnarray}
Further,
\begin{eqnarray*}
\lefteqn{{\bf P}\{A(y,z)\}}\\
&=& \overline B(x+ya+za) \overline B(x+z(a-b'))\\
&& +\ {\bf P}\Bigl\{x+ya+z(a-b')<\sigma_1\le x+ya+za,
\ \sigma_2>x+z(a-b'),\\
&& \hspace{60mm}
\sigma_1+\sigma_2> 2x+ya+z(2a-b')\Bigr\}\\
&=& \overline B(x+ya+za) \overline B(x+z(a-b'))\\
&& +\ {\bf P}\Bigl\{x+ya+z(a{-}b')<\sigma_1\le x+ya+za,\
\sigma_1+\sigma_2> 2x+ya+z(2a{-}b')\Bigr\}\\
&\equiv& \overline B(x+ya+za) \overline B(x+z(a-b'))+Q(y,z),
\end{eqnarray*}
since the event
$\{\sigma_1\le x+ya+za,\sigma_1+\sigma_2>2x+ya+z(2a-b')\}$
implies the event $\{\sigma_2>x+z(a'-b)\}$.
Consequently integrating over $y$ and $z$, we obtain
\begin{eqnarray*}
\int_0^\infty\int_0^\infty
\overline B(x{+}ya{+}za) \overline B(x{+}z(a{-}b')) dydz
&=& \frac{1}{a}
\int_0^\infty \overline B_I (x{+}za) \overline B(x{+}z(a{-}b'))dz.
\end{eqnarray*}
By the total probability formula,
\begin{eqnarray*}
Q(y,z) &=& \int_0^{zb'}
{\bf P}\{\sigma_1\in x+ya+z(a{-}b')+dt\}
{\bf P}\{\sigma_2>x+za-t\}\\
&=& \int_0^{zb'} \overline B(x+za-t) B(x+ya+z(a{-}b')+dt).
\end{eqnarray*}
The integration against $y$ leads to the equalities
\begin{eqnarray*}
\int_0^\infty Q(y,z)dy
&=& \frac{1}{a} \int_0^{zb'}
\overline B(x+za-t) B_I(x+z(a-b')+dt)\\
&=& \frac{1}{a} \int_0^{zb'}
\overline B(x+za-t) \overline B(x+z(a-b')+t)dt\\
&=& \frac{b'}{a} \int_0^z
\overline B(x+za-tb') \overline B(x+z(a-b')+tb')dt.
\end{eqnarray*}
Integrating against $z$, we obtain:
\begin{eqnarray*}
\int_0^\infty \int_0^\infty Q(y,z) dydz
&=& \frac{b'}{a} \int_0^\infty\int_0^z
\overline B(x+za-tb') \overline B(x+z(a-b')+tb')dtdz\\
&=& \frac{b'}{a} \int_0^\infty \int_t^\infty
\overline B(x+za-tb') \overline B(x+z(a-b')+tb')dzdt\\
&=& \frac{b'}{a} \int_0^\infty \int_0^\infty
\overline B(x+za+t(a{-}b')) \overline B(x+z(a{-}b')+ta)dzdt.
\end{eqnarray*}
By Lemma \ref{calcul.1} with $f(y)=\overline B(x+y)$,
$\alpha=a$, and $\beta=a-b'$, the latter integral is equal to
\begin{eqnarray*}
\frac{1}{a(2a-b')} (\overline B_I (x))^2
-\frac{2(a-b')}{a(2a-b')} \int_0^\infty
\overline B_I (x+ya) \overline B(x+y(a-b'))dy.
\end{eqnarray*}
Putting everything together into (\ref{integral.repr}),
we obtain the following equivalence, as $x\to\infty$:
\begin{eqnarray*}
\sum_{l=1}^\infty \sum_{k=1}^\infty {\bf P}\{A_{kl}'\}
&\sim& \frac{1}{a(2a-b')} (\overline B_I(x))^2\\
&& \hspace{10mm} +\frac{b'}{a(2a-b')} \int_0^\infty \overline B_I(x+ya)
\overline B(x+y(a-b'))dy,
\end{eqnarray*}
which due to (\ref{throw.A.prime}) completes the proof of Lemma.

{\bf \ref{sec.2.max.lower}.3. Proof of Theorem \ref{th.2.max.lower}.}
If $\overline B_I (x)$ is long-tailed,
then the function in $x$
\begin{eqnarray*}
(\overline B_I (x))^2+
b\int_0^\infty \overline B_I (x+ya)\overline B(x+y(a-b))dy
\end{eqnarray*}
is long-tailed as well. Indeed, for any fixed $t$,
we have, as $x\to\infty$,
\begin{eqnarray*}
\int_0^\infty \overline B_I (x{+}t{+}ya)
\overline B(x{+}t{+}y(a{-}b))dy
&\sim& \int_0^\infty \overline B_I (x{+}ya)
\overline B(x{+}t{+}y(a{-}b))dy.
\end{eqnarray*}
Integrating by parts we get the equality for RHS integral
\begin{eqnarray*}
\lefteqn{-\frac{1}{a-b}\overline B_I (x{+}ya)
\overline B_I(x{+}t{+}y(a{-}b))\Bigr|_0^\infty
-\int_0^\infty \overline B(x{+}ya)
\overline B_I (x{+}t{+}y(a{-}b))dy}\\
&& \hspace{30mm} \sim \frac{1}{a-b}(\overline B_I (x))^2
-\int_0^\infty \overline B(x{+}ya)
\overline B_I (x{+}y(a{-}b))dy\\
&& \hspace{60mm} = \int_0^\infty \overline B_I(x{+}ya)
\overline B (x{+}y(a{-}b))dy.
\end{eqnarray*}
So, we can apply Lemma \ref{lower.bound},
and it is sufficient to prove the lower bound of Theorem
\ref{th.2.max.lower} for the queueing system $D/GI/2$ with
deterministic input stream.
Let the interarrival times $\tau_n$ be deterministic,
i.e., $\tau_n\equiv a$.
Then the event $A_{nkl}$ implies the event
\begin{eqnarray*}
\lefteqn{\Bigl\{W_{k+1,2}>x+(l{-}k)a+(n{-}l)(a{-}b')-a,\
W_{l+1,1}>x+(n{-}l)(a{-}b')-a,}\\
&& \hspace{35mm}
W_{l+1,1}+W_{k+1,2}>2x+(l{-}k)a+(n{-}l)(2a{-}b')-2a\Bigr\},
\end{eqnarray*}
which implies
\begin{eqnarray*}
\Bigl\{W_{l+1,2}, W_{l+1,1}>x+(n{-}l)(a{-}b'){-}a,
\ W_{l+1,1}{+}W_{l+1,2}>2x+(n{-}l)(2a{-}b'){-}2a\Bigr\}.
\end{eqnarray*}
Thus, by Corollary \ref{SLLN.max.2}
(with $\xi=\sigma-\tau$ and $\eta=-\tau$),
there exists $N$ such that
\begin{eqnarray}\label{corr.to.corr2}
{\bf P}\{W_n>x \mid A_{nkl}\} &\ge& 1-\varepsilon
\end{eqnarray}
for any $n>N$ and $k<l<n-N$.

Taking into account that the events $A_{nkl}\cap C_{nkl}$
are disjoint for distinct pairs $(k,l)$,
we obtain the following estimates:
\begin{eqnarray*}
{\bf P}\{W_n>x\}
&\ge& \sum_{k=1}^{n-N} \sum_{l=k+1}^{n-N}
{\bf P}\{W_n>x,A_{nkl},C_{nkl}\}\\
&\ge& \sum_{k=1}^{n-N} \sum_{l=k+1}^{n-N}
{\bf P}\{W_n>x,A_{nkl}\}
-\sum_{k=1}^{n-N} \sum_{l=k+1}^{n-N}
{\bf P}\{A_{nkl},\overline C_{nkl}\}.
\end{eqnarray*}
Since the events $A_{nkl}$ and $C_{nkl}$ are independent,
\begin{eqnarray*}
{\bf P}\{W_n>x\}
&\ge& \sum_{k=1}^{n-N} \sum_{l=k+1}^{n-N}
{\bf P}\{W_n>x,A_{nkl}\}-\sup_{kl} {\bf P}\{C_{nkl}\}
\sum_{k=1}^{n-N} \sum_{l=k+1}^{n-N} {\bf P}\{A_{nkl}\}\\
&=& \sum_{k=1}^{n-N} \sum_{l=k+1}^{n-N}
{\bf P}\{W_n>x \mid A_{nkl}\}{\bf P}\{A_{nkl}\}
-o(1)\sum_{k=1}^{n-N} \sum_{l=k+1}^{n-N} {\bf P}\{A_{nkl}\}
\end{eqnarray*}
as $x\to\infty$ uniformly in $n$, by (\ref{slln.for.C}).
Together with (\ref{corr.to.corr2}) it implies that,
for sufficiently large $x$ and $n>N$,
\begin{eqnarray*}
{\bf P}\{W_n>x\} &\ge&
(1-2\varepsilon)\sum_{k=1}^{n-N}
\sum_{l=k+1}^{n-N} {\bf P}\{A_{nkl}\}.
\end{eqnarray*}
Letting now $n\to\infty$, we derive from
Lemma \ref{int.2} the following lower bound,
for all sufficiently large $x$:
\begin{eqnarray*}
{\bf P}\{W>x\}
&\ge& \frac{1-3\varepsilon}{a(2a-b')} \Bigl[(\overline B_I(x))^2
+ b'\int_0^\infty \overline B_I(x+ya)
\overline B(x+y(a-b'))dy\Bigr].
\end{eqnarray*}
Note that, for any $b'<b<a$,
$$
\int_0^\infty \overline B_I(x+ya) \overline B(x+y(a-b'))dy
\ge \frac{a-b}{a-b'}
\int_0^\infty \overline B_I(x+ya)
\overline B(x+y(a-b))dy.
$$
We complete the proof of the Theorem by letting
$\varepsilon\downarrow 0$.

\mysection{The maximal stability case: an upper bound
\label{sec.2.max.upper}}

\begin{Theorem}\label{th.2.max.upper}
Assume $b\in (0,a)$.
Suppose that the distribution $B_I$ is subexponential.
Then, as $x\to\infty$,
\begin{eqnarray*}
{\bf P}\{W>x\}
&\le& \frac{1{+}o(1)}{a(2a{-}b)} \Biggl[(\overline B_I(x))^2
+ b \int_0^\infty \overline B_I(x{+}ya)
\overline B(x{+}y(a{-}b))dy\Biggr].
\end{eqnarray*}
\end{Theorem}

By Lemma \ref{upper.bound}, it is sufficient to prove this upper
bound for the queueing system $D/GI/2$ with deterministic input stream.
So, let the interarrival times $\tau_n$ be deterministic,
i.e., $\tau_n\equiv a$.
Let $\sigma_n^{(1)}$ and $\sigma_n^{(2)}$, $n\ge1$,
be independent random variables with common distribution $B$.
In this Section, define the service times $\sigma_n$ recursively.
For that, we have to associate workloads with servers. Put
$U_1=(U_{1,1},U_{1,2})=(0,0)$ and introduce the recursion
\begin{eqnarray}\label{UU}
U_{n+1} = (U_n + e_{\alpha_n}\sigma_n - i a)^+
\end{eqnarray}
where $\alpha_n=1$ if $U_{n,1}<U_{n,2}$ and $\alpha_n=2$ if
$U_{n,1}>U_{n,2}$. If $U_{n,1}=U_{n,2}$, then $\alpha_n$ takes
values $1$ and $2$ with equal probabilities independently of
everything else.
Note that  $W_n = R(U_n)$ a.s. for any $n=1$, 2, \ldots.

Now we can define $\sigma_n$ by induction. Indeed,
$\alpha_0$ is chosen at random from the set $\{1 , 2\}$.
Put $\sigma_0 =
\sigma_0^{(\alpha_0)}$. Then $U_1$ is defined by recursion
(\ref{UU}) with $n=0$.  Assume that $U_n$ is defined for some
$n>0$. Then $\alpha_n$ is defined, too. Put $\sigma_n=
\sigma_n^{(\alpha_n)}$ and determine $U_{n+1}$ by (\ref{UU}).

Due to the symmetry, for any $n$,
\begin{eqnarray}\label{symmetry}
{\bf P}\{\alpha_n=1\} &=&
{\bf P}\{\alpha_n=2 \}=1/2.
\end{eqnarray}

Consider two auxiliary
$D/GI/1$ queueing systems which work in parallel:
at any time instant $T_n = na$, $n=1$, 2, \ldots,
one customer arrives in the first queue and one in the second.
Service times in queue $i=1,2$ are equal to $\sigma_n^{(i)}$.
Denote by $W_n^{(i)}$, $i=1$, $2$, the waiting
times in the $i$th queue and put $W_n^{(1)}=W_n^{(2)}=0$.
Since $b<a$, both queues are stable.
Let $W^{(i)}$ be a stationary waiting time in the $i$th queue.
By monotonicity, with probability 1,
\begin{eqnarray}\label{W.le.min.W1.W2.prest}
W_n &\le& \min\,(W_n^{(1)},\,W_n^{(2)})
\end{eqnarray}
for any $n\ge1$. Hence,
\begin{eqnarray}\label{W.le.min.W1.W2}
W &\le& \min\,\{W^{(1)},\,W^{(2)}\}.
\end{eqnarray}

\begin{Lemma}\label{max_upper_independence}
The waiting times $\{W_n^{(1)}\}$ and $\{W_n^{(2)}\}$ are independent.
\end{Lemma}

\proof\ \ follows from the observation that the input
(deterministic) stream and service times in
the first queue do not depend on the input
(also deterministic) stream and
service times in the second one.

Provided $B_I$ is a subexponential distribution,
\begin{eqnarray}\label{W1.sim.BI}
{\bf P}\{W^{(i)}>x\} &\sim& \frac{1}{a-b}\overline B_I(x)
\quad\mbox{ as } x\to\infty.
\end{eqnarray}
Then Lemma \ref{max_upper_independence} together with
(\ref{W.le.min.W1.W2}) implies the following
simple upper bound:
\begin{eqnarray}\label{W.le.W1sq}
\limsup_{x\to\infty} \frac{{\bf P}\{W>x\}}{(\overline B_I(x))^2}
&\le& \frac{1}{(a-b)^2}.
\end{eqnarray}

\remark\label{rem.s.min.up}
For a $GI/GI/s$ queue with $a<b$ and subexponential
distribution $B_I$, similar arguments lead to
$$
\limsup_{x\to\infty}
\frac{{\bf P}\{W>x\}}{(\overline{B}_I(x))^s}
\le \frac{1}{(a-b)^s}.
$$

Introduce the events, for $k<n$,
\begin{eqnarray*}
A_{nk}^{(1)} &=& \{\sigma_k^{(1)}>x+(n-k)(a-b)\},\\
A_{nk}^{(2)} &=& \{\sigma_k^{(2)}>x+(n-k)(a-b)\}.
\end{eqnarray*}

\begin{Lemma}[See also {[}\ref{BaF}, Theorem 5{]}]
\label{upper_1server}
Provided the distribution $B_I$ is subexponential,
for any fixed $N$,
\begin{eqnarray*}
\limsup_{n\to\infty}
{\bf P}\Bigl\{W^{(1)}_n>x,\,
\bigcap_{k=1}^{n-N}\overline{A^{(1)}_{nk}}\Bigr\}
&=& o(\overline B_I (x))
\quad\mbox{ as }x\to\infty.
\end{eqnarray*}
\end{Lemma}

\proof.
For any $\delta>0$, consider the disjoint events
\begin{eqnarray*}
C_{nk}^{(1)} &=&
\Bigl\{\Bigl\{\sigma_k^{(1)}>x+(n-k)(a-b+\delta)\Bigr\}\cap
\bigcap_{\stackrel{j=1}{j\neq k}}^{n-1}\Bigl\{\sigma_j^{(1)}\le x+(n-j)(a-b)\Bigr\}\Bigr\}.
\end{eqnarray*}
Due to the Law of Large Numbers, there exists $M>N$ such that
\begin{eqnarray*}
{\bf P}\{W^{(1)}_n>x\mid C_{nk}^{(1)}\} &\ge& 1-\delta
\end{eqnarray*}
for any $n\ge M$ and $k\le n-M$ and, by the limit at (\ref{slln.for.C}),
\begin{eqnarray*}
{\bf P}\{C_{nk}^{(1)}\}
&\ge& (1-\delta){\bf P}\{\sigma_k^{(1)}>x+(n-k)(a-b+\delta)\}.
\end{eqnarray*}
The events $C_{nk}^{(1)}$, $k\le n-M$, are disjoint, hence,
\begin{eqnarray*}
{\bf P}\Bigl\{W^{(1)}_n>x,\,\bigcup_{k=1}^{n-M} C^{(1)}_{nk}\Bigr\}
&=& \sum_{k=1}^{n-M} {\bf P}\{W^{(1)}_n>x,\,C^{(1)}_{nk}\}\\
&\ge& (1-\delta)^2\sum_{k=M}^{n-1}
{\bf P}\{\sigma_k^{(1)}>x+k(a-b+\delta))\}.
\end{eqnarray*}
The latter implies the following lower bound:
\begin{eqnarray*}
\liminf_{n\to\infty}
{\bf P}\Bigl\{W^{(1)}_n>x,\,\bigcup_{k=1}^{n-M} C^{(1)}_{nk}\Bigr\}
&\ge& (1-\delta)^2\sum_{k=M}^\infty
\overline B(x+k(a-b+\delta))\\
&\sim& \frac{(1-\delta)^2}{a-b+\delta}\overline B_I(x)
\end{eqnarray*}
as $x\to\infty$. Since $A^{(1)}_{nk}\supseteq C^{(1)}_{nk}$
and since
$M>N$ and $\delta>0$ can be chosen arbitrarily,
\begin{eqnarray*}
\liminf_{n\to\infty}
{\bf P}\Bigl\{W^{(1)}_n>x,\,\bigcup_{k=1}^{n-N} A^{(1)}_{nk}\Bigr\}
&\ge& \frac{1+o(1)}{a-b}\overline B_I(x)
\quad\mbox{ as }x\to\infty.
\end{eqnarray*}
Together with (\ref{W1.sim.BI}), it implies the assertion of Lemma.

\proof\ of Theorem \ref{th.2.max.upper} continued.
Estimate (\ref{W.le.min.W1.W2.prest})
and Lemma \ref{max_upper_independence} imply
\begin{eqnarray*}
{\bf P}\Bigl\{W_n>x,\,
\bigcap_{k=1}^{n-N}\overline{A^{(1)}_{nk}}\cup
\bigcap_{l=1}^{n-N}\overline{A^{(2)}_{nl}}\Bigr\}
&\le& {\bf P}\Bigl\{W^{(1)}_n>x,\,W^{(2)}_n>x,\,
\bigcap_{k=1}^{n-N}\overline{A^{(1)}_{nk}}\cup
\bigcap_{l=1}^{n-N}\overline{A^{(2)}_{nl}}\Bigr\}\\
&\le& {\bf P}\Bigl\{W^{(1)}_n>x,\,
\bigcap_{k=1}^{n-N}\overline{A^{(1)}_{nk}}\Bigr\}
{\bf P}\{W^{(2)}_n>x\}\\
&& + {\bf P}\{W^{(1)}_n>x\}
{\bf P}\Bigl\{W^{(2)}_n>x,\,
\bigcap_{l=1}^{n-N}\overline{A^{(2)}_{nl}}\Bigr\}.
\end{eqnarray*}
Applying now Lemma \ref{upper_1server} and relation
(\ref{W1.sim.BI}), we conclude that,
as $x\to\infty$,
\begin{eqnarray*}
\limsup_{n\to\infty}
{\bf P}\Bigl\{W_n>x,\,
\bigcap_{k=1}^{n-N}\overline{A_{nk}^{(1)}}\cup
\bigcap_{l=1}^{n-N}\overline{A_{nl}^{(2)}}\Bigr\}
&=& o((\overline B_I(x))^2).
\end{eqnarray*}
Since
\begin{eqnarray*}
\bigcap_{k=1}^{n-N}\overline{A_{nk}^{(1)}}\cup
\bigcap_{l=1}^{n-N}\overline{A_{nl}^{(2)}}
&=& \bigcap_{k,l=1}^{n-N}\Bigl(\overline{A_{nk}^{(1)}}\cup
\overline{A_{nl}^{(2)}}\Bigr)
= \overline{\bigcup_{k,l=1}^{n-N}\Bigl(A_{nk}^{(1)}\cap
A_{nl}^{(2)}\Bigr)},
\end{eqnarray*}
we obtain the equivalent relation, as $x\to\infty$,
\begin{eqnarray}\label{upper_2server}
\limsup_{n\to\infty}
{\bf P}\Biggl\{W_n>x,\,
\overline{\bigcup_{k,l=1}^{n-N}\Bigl(A_{nk}^{(1)}\cap
A_{nl}^{(2)}\Bigr)}\Biggr\}
&=& o((\overline B_I(x))^2).
\end{eqnarray}

Fix $\varepsilon>0$ and put $b'=b+\varepsilon$.
For any $n$ and $k\le l\le n$, define
\begin{eqnarray*}
D^{(1)}_{nk} &=& \{\sigma^{(1)}_k>x+(l-k)a+(n-l)(a-b')\},\\
D^{(2)}_{nl} &=& \{\sigma^{(2)}_l>x+(n-l)(a-b')\},\\
D_{nkl} &=&
\{\sigma^{(1)}_k+\sigma^{(2)}_l>2x+(l-k)a+(n-l)(2a-b')\}.
\end{eqnarray*}
For any $n$ and $l\le k\le n$, define
\begin{eqnarray*}
D^{(1)}_{nk} &=& \{\sigma^{(1)}_k>x+(n-k)(a-b')\},\\
D^{(2)}_{nl} &=& \{\sigma^{(2)}_l>x+(k-l)a+(n-k)(a-b')\},\\
D_{nkl} &=&
\{\sigma^{(1)}_k+\sigma^{(2)}_l>2x+(k-l)a+(n-k)(2a-b')\}.
\end{eqnarray*}
Denote
\begin{eqnarray*}
F_{nkl} &=& D^{(1)}_{nk}\cap D^{(2)}_{nl}\cap D_{nkl}.
\end{eqnarray*}
We can derive an upper bound on the probability of the event
$\{W_n>x\}$ as follows:
\begin{eqnarray}\label{P1.P2.P3.dec}
\lefteqn{{\bf P}\{W_n>x\}}\nonumber\\
&\le& {\bf P}\Bigl\{W_n>x,\,\bigcup_{k,l=1}^{n-N} F_{nkl}\Bigr\}
+ {\bf P}\Bigl\{W_n>x,\,
\overline{\bigcup_{k,l=1}^{n-N} F_{nkl}},\,
\bigcup_{k,l=1}^{n-N}\Bigl(A_{nk}^{(1)}\cap
A_{nl}^{(2)}\Bigr)\Bigr\}\nonumber\\
&& + {\bf P}\Bigl\{W_n>x,\,
\overline{\bigcup_{k,l=1}^{n-N}\Bigl(A_{nk}^{(1)}\cap
A_{nl}^{(2)}\Bigr)}\Bigr\}\nonumber\\
&\equiv& P_{n1}+P_{n2}+P_{n3}.
\end{eqnarray}
Here the first term is not greater than
\begin{eqnarray}\label{dec.of.p1}
P_{n1} &\le&
{\bf P}\Bigl\{W_n>x,\,
\bigcup_{\stackrel{k,l=1}{k<l}}^{n-1} F_{nkl}\Bigr\}
+ {\bf P}\Bigl\{W_n>x,\,
\bigcup_{\stackrel{k,l=1}{k>l}}^{n-1} F_{nkl}\Bigr\}
+ {\bf P}\Bigl\{W_n>x,\,
\bigcup_{k=1}^{n-1}F_{nkk}\Bigr\}\nonumber\\
&\equiv& P_{n11}+P_{n12}+P_{n13}.
\end{eqnarray}
The third probability is negligible in the sense that
\begin{eqnarray}\label{p13}
P_{n13} &\le& {\bf P}\Bigl\{\bigcup_{k=1}^{n-1}
(D^{(1)}_{nk}\cap D^{(2)}_{nk})\Bigr\}
\le\sum_{k=1}^{n-1} {\bf P}\{D^{(1)}_{nk}\}
{\bf P}\{D^{(2)}_{nk}\}\nonumber\\
&\le& \overline B(x)\sum_{k=1}^\infty
\overline B(x+k(a-b-\varepsilon))\nonumber\\
&\le& \overline B(x)\overline B_I(x)/(a-b-\varepsilon)
= o((\overline B_I(x))^2)
\end{eqnarray}
as $x\to\infty$, since $\overline B(x)=o(\overline B_I (x))$.
The first probability in (\ref{dec.of.p1}) admits the
following upper bound:
\begin{eqnarray*}
P_{n11} &\le& \sum_{k=1}^{n-1}
{\bf P}\Bigl\{W_n>x,\,D^{(1)}_{nk}, \alpha_k=1,
\bigcup_{l=k+1}^{n-1}(D^{(2)}_{nl}\cap D_{nkl})\Bigr\}\\
&& + \sum_{k=1}^{n-1}
{\bf P}\Bigl\{W_n>x,\,D^{(1)}_{nk}, \alpha_k=2,
\bigcup_{l=k+1}^{n-1}(D^{(2)}_{nl}\cap D_{nkl})\Bigr\}
\equiv\Sigma_1+\Sigma_2.
\end{eqnarray*}
For $\Sigma_1$, we have the following inequality
and equalities:
\begin{eqnarray}\label{sigma.1}
\Sigma_1 &\le& \sum_{\stackrel{k,l=1}{k<l}}^{n-1}
{\bf P}\Bigl\{D^{(1)}_{nk},
\alpha_k=1,D^{(2)}_{nl},D_{nkl}\Bigr\}\nonumber\\
&=& \sum_{\stackrel{k,l=1}{k<l}}^{n-1}
{\bf P}\{\alpha_k=1\}
{\bf P}\Bigl\{D^{(1)}_{nk},D^{(2)}_{nl},D_{nkl}\Bigr\}
= \frac{1}{2}\sum_{\stackrel{k,l=1}{k<l}}^{n-1}
{\bf P}\{F_{nkl}\},
\end{eqnarray}
by independence of the event $\{ \alpha_k=1 \}$
from $D^{(1)}_{nk}$, $D^{(2)}_{nl}$ and $D_{nkl}$
and by the symmetry (\ref{symmetry}).
The sum $\Sigma_2$ is not greater than
\begin{eqnarray*}
\Sigma_2 &\le& \sum_{k=1}^{n-1}
{\bf P}\Bigl\{W_n>x,\,D^{(1)}_{nk},
\alpha_k=2\Bigr\}\\
&=& \sum_{k=1}^{n-1} {\bf P}\{D^{(1)}_{nk}\}
{\bf P}\Bigl\{W_n>x,\,
\alpha_k=2\Bigr\}\le {\bf P}\{W_n>x\}
\sum_{k=1}^{n-1} {\bf P}\{D^{(1)}_{nk}\}.
\end{eqnarray*}
Hence, $\Sigma_2=o({\bf P}\{W_n>x\})$ as $x\to\infty$
uniformly in $n$. Combining the latter fact with estimate (\ref{sigma.1})
for $\Sigma_1$, we get
\begin{eqnarray}\label{p11}
P_{n11} &\le& \frac{1}{2}\sum_{\stackrel{k,l=1}{k<l}}^{n-1}
{\bf P}\{F_{nkl}\}+o({\bf P}\{W_n>x\}).
\end{eqnarray}
Taking into account the equality $P_{n11}=P_{n12}$,
we obtain from (\ref{dec.of.p1}), (\ref{p13}) and
(\ref{p11}) the following estimate:
\begin{eqnarray*}
P_{n1} &\le& \sum_{\stackrel{k,l=1}{k<l}}^{n-1}
{\bf P}\{F_{nkl}\}+o((\overline B_I(x))^2)
\end{eqnarray*}
as $x\to\infty$ uniformly in $n$.
Now applying the calculations of Section 3.3
we can write down the following estimate, as $x\to\infty$:
\begin{eqnarray}\label{P1.dec}
\limsup_{n\to\infty} P_{n1}
&\le& \frac{1+o(1)}{a(2a-b')}\Bigl[(\overline B_I (x))^2
+b'\int_0^\infty\overline B_I(x+ya')\overline B(x+y(a-b'))dy\Bigr].
\nonumber\\
\end{eqnarray}

It is proved in (\ref{upper_2server}) that, uniformly in $n$,
\begin{eqnarray}\label{P2.dec}
P_{n3} &=& o((\overline B_I(x))^2) \quad\mbox{ as }x\to\infty.
\end{eqnarray}

We have
\begin{eqnarray*}
\overline{\bigcup_{k,l=1}^{n-N} F_{nkl}} \cap
\bigcup_{k,l=1}^{n-N}\Bigl(A_{nk}^{(1)}\cap A_{nl}^{(2)}\Bigr)
&\subseteq& \bigcup_{k,l=1}^{n-N}\Bigl(A_{nk}^{(1)}\cap
A_{nl}^{(2)}\cap\overline F_{nkl}\Bigr).
\end{eqnarray*}
Thus,
\begin{eqnarray}\label{P3.dec.pr}
P_{n2} &\le& \sum_{k,l=1}^{n-N}
{\bf P}\{W_n>x \mid A_{nk}^{(1)},A_{nl}^{(2)},\overline F_{nkl}\}
{\bf P}\{A_{nk}^{(1)}\cap A_{nl}^{(2)}\}.
\end{eqnarray}
Conditioning on $W_{nk}$ and $W_{nl}$ yields,
for any $w>0$,
\begin{eqnarray*}
{\bf P}\{W_n>x \mid A_{nk}^{(1)},A_{nl}^{(2)},\overline F_{nkl}\}
&\le& {\bf P}\{W_n>x \mid W_{k1}\le w,W_{l2}\le w,
A_{nk}^{(1)},A_{nl}^{(2)},\overline F_{nkl}\}\\
&& + {\bf P}\{W_{k1}>w\}+{\bf P}\{W_{l2}>w\}.
\end{eqnarray*}
Since $b<2a$, the two-server queue is stable and, in particular, the
sequence of distributions of random variables $(W_{n1},W_{n2})$ is
tight.
It means that,
for any fixed $\varepsilon>0$, there exists $w$ such that,
for any $k\ge0$ and $l\ge0$,
\begin{eqnarray*}
{\bf P}\{W_{k1}>w\} &\le& \varepsilon \quad
\mbox{ and }\quad {\bf P}\{W_{l2}>w\} \le \varepsilon.
\end{eqnarray*}
Also, from the stability and from
Corollary \ref{SLLN.max.3},
for any fixed $\varepsilon>0$ and $w>0$,
there exists $N$ such that,
for any $n\ge N$ and $k$, $l\le n-N$,
\begin{eqnarray*}
{\bf P}\{W_n>x \mid W_{k1}\le w,W_{l2}\le w,
A_{nk}^{(1)},A_{nl}^{(2)},\overline F_{nkl}\}
&\le& \varepsilon.
\end{eqnarray*}
Combining these estimates we obtain from (\ref{P3.dec.pr}),
\begin{eqnarray*}
P_{n2} &\le& 3\varepsilon \sum_{k,l=1}^{n-N}
{\bf P}\{A_{nk}^{(1)}\cap A_{nl}^{(2)}\}
=3\varepsilon \Bigl(\sum_{k=1}^{n-N}
{\bf P}\{A_{nk}^{(1)}\}\Bigr)^2.
\end{eqnarray*}
Hence,
\begin{eqnarray}\label{P3.dec}
P_{n2} &\le& 3\varepsilon \Bigl(\sum_{k=1}^\infty
\overline B(x+k(a-b'))\Bigr)^2
\le \frac{3\varepsilon}{(a-b')^2}(\overline B_I(x))^2.
\end{eqnarray}
Since the choice of $\varepsilon>0$ is arbitrary,
relations (\ref{P1.P2.P3.dec})--(\ref{P2.dec})
and (\ref{P3.dec}) imply the conclusion of Theorem
\ref{th.2.max.upper}.

\mysection{The minimal stability case:
lower bounds\label{sec.2.min.lower}}

\begin{Theorem}\label{th.2.min.lower}
Let $b\in(a,2a)$ and the integrated tail distribution
$B_I$ be long tailed.
Then the tail of the stationary waiting time
satisfies the following inequality,
for any fixed $\delta>0$:
\begin{eqnarray*}
{\bf P}\{W>x\}
&\ge& \frac{1+o(1)}{2a-b}
\overline B_I\Bigl(\frac{b+\delta}{b-a}x\Bigr)
\quad \mbox{ as }x\to\infty.
\end{eqnarray*}
\end{Theorem}

Notice that if $b\in(a,2a)$
then $\frac{b}{b-a}>2$.

\remark\label{sserverlower}
By use of similar arguments, one can get the following
result for an $s$-server queue, $s\ge2$:
if the integrated distribution $B_I$ is long tailed and
$b\in((s-1)a,sa)$, then, for any $\delta>0$,
\begin{eqnarray*}
{\bf P}\{W>x\} &\ge&
\frac{1+o(1)}{sa-b} \overline B_I
\left(\frac{(s-1)b-s(s-2)a+\delta}{b-(s-1)a}x\right)
\quad \mbox{ as }x\to\infty.
\end{eqnarray*}

Theorem \ref{th.2.min.lower} implies the following
\begin{Corollary}\label{min_lower.reg}
Assume that $B_I \in {\cal IRV}$. Then, as $x\to\infty$,
\begin{eqnarray*}
{\bf P}\{W>x\} &\ge& \frac{1+o(1)}{2a-b}
\overline B_I\Bigl(\frac{b}{b-a}x\Bigr).
\end{eqnarray*}
\end{Corollary}

In the case $b\in[a,2a)$, one can also derive a lower bound
which is similar to (\ref{firstlb}).
More precisely, assume $b\in[a,2a)$.
Then introduce another two-server queue with the same
service times and with inter-arrival times
$\tau'_n=c\tau_n$, where $c>b/a$.
For this queue, denote by $W'$ a stationary waiting time
of a typical customer. Due to monotonicity,
${\bf P}\{W'>x\} \le {\bf P}\{W>x\}$ for all $x$.
Applying Theorem \ref{th.2.max.lower} and
Remark \ref{rem.rou}, we get the following lower bound
for the case $b\in[a,2a)$: if the integrated tail
distribution $B_I$ is long-tailed, then, for any $c>b/a$,
\begin{eqnarray}\label{anotherlower}
{\bf P}\{W>x\} &\ge&
(1+o(1))\frac{2ca+b}{2c^2a^2(2ca-b)}(\overline B_I (x))^2.
\end{eqnarray}

\proof\ of Theorem \ref{th.2.min.lower}.
By Lemma \ref{lower.bound}, it is sufficient to prove
the lower bound for the queueing system $D/GI/2$ with
deterministic input stream.
Let the interarrival times $\tau_n$ be deterministic,
i.e., $\tau_n\equiv a$. For any $\delta>0$, set
$\varepsilon=\frac{\delta(b-a)}{a+\delta}$.
Put $b'=b-\varepsilon$ and $N=\frac{x}{b'-a}$.
For any $k\in[1,n-N]$, consider the events
\begin{eqnarray*}
A_{nk} &=& \{\sigma_k>2x+(2a-b')(n-k)\},\\
C_{nk} &=& \bigcap_{\stackrel{l=1}{l\ne k}}^n
\{\sigma_l\le 2x+(2a-b')(n-l)\}.
\end{eqnarray*}
Since ${\bf E}\sigma$ is finite,
\begin{eqnarray}\label{c.to.1}
{\bf P}\{\overline C_{nk}\}
&\le& \sum_{l=1}^\infty
{\bf P}\{\sigma_1>2x+(2a-b')l\}
= O(\overline B_I (2x)) \to 0
\end{eqnarray}
as $x\to\infty$ uniformly in $n\ge1$ and $k\le n$.
Since the events $A_{nk}\cap C_{nk}$, $k\in[1,n]$,
are disjoint, we obtain
\begin{eqnarray*}
{\bf P}\{W_n>x\}
&\ge& \sum_{k=1}^{n-N}{\bf P}\{W_n>x,A_{nk},C_{nk}\}\\
&\ge& \sum_{k=1}^{n-N}{\bf P}\{W_n>x,A_{nk}\}-
\sum_{k=1}^{n-N}{\bf P}\{A_{nk},\overline C_{nk}\}.
\end{eqnarray*}
Since the events $A_{nk}$ and $C_{nk}$ are independent,
we get
\begin{eqnarray}\label{pequ}
{\bf P}\{W_n>x\}
&\ge& \sum_{k=1}^{n-N}{\bf P}\{W_n>x,A_{nk}\}-
\sup_{k\le n}{\bf P}\{\overline C_{nk}\}
\sum_{k=1}^{n-N}{\bf P}\{A_{nk}\}\nonumber\\
&=& \sum_{k=1}^{n-N}{\bf P}\{W_n>x\mid A_{nk}\}
{\bf P}\{A_{nk}\}-
o(1)\sum_{k=1}^{n-N}{\bf P}\{A_{nk}\}
\end{eqnarray}
as $x\to\infty$ uniformly in $n\ge1$, by (\ref{c.to.1}).
The event $A_{nk}$ implies the event
$$
W_{k+1,2}>2x+(2a-b')(n-k)-a.
$$
Thus, it follows from Corollary \ref{SLLN.max.4} that
\begin{eqnarray*}
{\bf P}\{W_n>x \mid A_{nk}\} &\to& 1
\end{eqnarray*}
as $x\to\infty$ uniformly in $n$ and $k\le n-N$.
Therefore, we can derive from (\ref{pequ}) the estimate
\begin{eqnarray*}
{\bf P}\{W>x\}=\lim_{n\to\infty} {\bf P}\{W_n>x\}
&\ge& (1-\varepsilon)\lim_{n\to\infty}
\sum_{k=1}^{n-N}{\bf P}\{A_{nk}\}\\
&=& (1-\varepsilon)\sum_{k=N}^\infty
\overline B(2x+(2a-b')k),
\end{eqnarray*}
which is valid for all sufficiently large $x$.
Since the tail $\overline B_I(v)$ is long-tailed,
\begin{eqnarray*}
\sum_{k=N}^\infty \overline B(2x+(2a-b')k)
&\sim& \frac{1}{2a-b'} \overline B_I(2x+(2a-b')N)\\
&=& \frac{1}{2a-b'} \overline B_I \Bigl(\frac{b'}{b'-a}x\Bigr)
= \frac{1}{2a-b'} \overline B_I \Bigl(\frac{b+\delta}{b-a}x\Bigr)
\end{eqnarray*}
as $x\to\infty$. The proof is complete.

\mysection{The minimal stability case: an upper bound}
\label{sec.2.min.upper}

\begin{Theorem} \label{min_upper}
Assume $b\in[a,2a)$. Let both
$B$ and $B_I$ be subexponential distributions.
Then the tail of the stationary waiting time
satisfies the following inequality,
as $x\to\infty${\rm:}
\begin{eqnarray*}
{\bf P}\{W\ge x\}
&\le& \frac{1+o(1)}{2a-b} \overline{B}_I(2x).
\end{eqnarray*}
\end{Theorem}

\remark\label{sserverupper}
By use of the same arguments, one can get the following
result for any $s$-server queue, $s\ge2$:
if $B_I\in{\mathcal S}$ and $b<sa$, then
\begin{eqnarray*}
{\bf P}\{W>x\} &\le& \frac{1+o(1)}{sa-b}
\overline{B}_I
(sx)
\quad \mbox{ as }x\to\infty
\end{eqnarray*}
provided that either (i)  $\sigma_1\ge (s-1)a$ a.s.,
or (ii) $B\in{\mathcal S}$.

\remark \label{FChupper}
For an $s$-server queue, Foss and Chernova [\ref{FCh}]
have proposed another way of obtaining upper bounds;
it is based on comparison with a queue with the so-called
{\it cyclic} service discipline.

\proof\ of Theorem \ref{min_upper}. From Lemma \ref{SLLN.max},
it is sufficient to consider the case of constant
interarrival times $\tau_n\equiv a$ only.
Put $M_{n,0}=0$ and
$$
M_{n,i+1} = (M_{n,i}+\sigma_{n+i}-a)^+.
$$
Since $b>a$, $M_{0,n}\to\infty$ a.s. as $n\to\infty$ and,
due to the Law of Large Numbers,
\begin{equation} \label{SSL}
\frac{M_{0,n}}{n} \to b-a\quad \mbox{a.s.}
\end{equation}
and in mean. Note that ${\bf E}M_{0,n}\ge n(b-a)$,
since $M_{0,n}\ge\sigma_0+\ldots+\sigma_{n-1}-na$.
For any given $\varepsilon>0$,
choose an integer $L>0$ such that
\begin{equation}\label{LL}
\frac{{\bf E} M_{0,L}}{L} \in [b-a, b-a+\varepsilon ).
\end{equation}

Consider any initial workload vector $W_0=(W_{0,1}, W_{0,2}) \ge 0$.
Put $Z_n = W_{n,1}+W_{n,2}$.
Since the increments of the minimal coordinate of
the waiting time vector is not greater than the
increments of $M_{0,n}$,
\begin{eqnarray*}
W_{1,n}-W_{1,0} &\le& M_{0,n}\quad\mbox{ for any }n.
\end{eqnarray*}
Hence, provided $W_{n,2}\ge a$, we have the inequality
\begin{eqnarray*}
Z_{n+1}-Z_n &\le& M_{0,n+1}-M_{0,n}-a.
\end{eqnarray*}
If $Z_0\ge 2aL$, then $W_{0,2}\ge aL$ and,
for $n=0$, \ldots, $L-1$, $W_{n,2}\ge a(L-n)\ge a$.
Therefore, if $Z_0 \ge 2aL$, then
\begin{eqnarray*}
Z_L &\le& Z_0+M_{0,L}-aL.
\end{eqnarray*}
Monotonicity implies, for any initial vector $W_0$,
\begin{eqnarray*}
Z_L &\le& \max \{2aL,Z_0\}+M_{0,L}-aL.
\end{eqnarray*}
Thus, the following inequalities are valid for any $n$:
\begin{eqnarray}\label{recu}
Z_{(n+1)L} &\le& \max \{2aL,Z_{nL}\}+M_{nL,L}-aL.
\end{eqnarray}

Consider a single-server queue with i.i.d.\
service times $\widehat{\sigma}_n= M_{nL,L}$ and constant
inter-arrival times $\widehat\tau_n=La$ and denote
by $\widehat W_n$ a waiting time of $n$th customer.
This queue is stable since
$\widehat b\equiv{\bf E}\widehat \sigma_1<aL\equiv\widehat a$.
Put $\widehat W_0=0$. Assuming that $Z_0=0$,
we can derive from (\ref{recu}) the following bounds:
for all $n=0$, 1, \ldots,
\begin{eqnarray}\label{Z.nL.W}
Z_{nL} &\le& 2aL+\widehat W_n \quad \mbox{a.s.}
\end{eqnarray}

Denote $\overline G(x)={\bf P}\{\widehat\sigma_0>x\}$.
We show that integrated tail distribution $G_I$ is
subexponential one. We need to consider only the case $L>1$.
Note first that
\begin{equation}\label{first}
\sigma_0+\ldots + \sigma_{L-1}-La \le \widehat\sigma_0
\leq \sigma_0+\ldots+\sigma_{L-1} \quad \mbox{a.s.}
\end{equation}
Since the distribution of $\sigma_1$ is assumed
to be subexponential, the asymptotics for the lower and
upper bounds in the latter inequalities are the same:
as $x\to\infty$,
\begin{equation}\label{second}
{\bf P}\Bigl\{\sum_0^{L-1}\sigma_i-La>x\Bigr\}
\sim {\bf P}\Bigl\{\sum_0^{L-1}\sigma_i>x\Bigr\}
\sim L\overline B(x).
\end{equation}
Therefore, the tail $\overline G(x)$ has the same
asymptotics and $G$ is a subexponential distribution. Thus,
\begin{eqnarray}\label{G_I.via.B_I}
\overline G_I(x) &=& \int_x^\infty \overline G(y)dy
\sim L\overline B_I(x).
\end{eqnarray}
and $G_I$ is a subexponential distribution, too.
Thus, by classic result (\ref{W.single}) for the single
server queue, the steady state distribution of the waiting
time $\widehat W_n$ satisfies the following relations,
as $x\to\infty$:
\begin{eqnarray}\label{asy.for.hat.W}
\lim_{n\to\infty}{\bf P}\{\widehat W_n>x\}
&\sim& \frac{1}{\widehat a{-}\widehat b}\overline G_I(x)
\le \frac{1}{(2a{-}b{-}\varepsilon)L}\overline G_I(x)
\sim \frac{1}{2a{-}b{-}\varepsilon}\overline B_I(x),
\end{eqnarray}
by (\ref{LL}) and (\ref{G_I.via.B_I}).
Since $Z_n=W_{n,1}+W_{n,2}\ge 2W_{n,1}$,
\begin{eqnarray*}
{\bf P}\{W>x\} &=& {\bf P}\{2W>2x\}
\le \lim_{n\to\infty}{\bf P}\{Z_{nL}>2x\}.
\end{eqnarray*}
Now it follows from (\ref{Z.nL.W}) and
(\ref{asy.for.hat.W}) that
\begin{eqnarray*}
{\bf P}\{W>x\}
&\le& \lim_{n\to\infty}{\bf P}\{\widehat W_n>2x-2aL\}\\
&\le& \frac{1+o(1)}{2a-b-\varepsilon}\overline B_I(2x-2aL)
\sim \frac{1}{2a-b-\varepsilon}\overline B_I(2x),
\end{eqnarray*}
since $B_I$ is long-tailed.
Letting $\varepsilon\downarrow 0$ concludes the proof.

\mysection{The minimal stability case: exact asymptotics}
\label{min.stab.proof}

In this Section, we prove Theorem \ref{th.2.min.reg}.
First note that, as follows from
(\ref{anotherlower}), the tail ${\bf  P}\{W>x\}$
may be heavier than that in Theorem  \ref{th.2.min.reg}, in general.
For instance, this happens if
\begin{equation}\label{ppp}
\overline B_I\Bigl(\frac{b}{b-a}x\Bigr)
=o(\overline B_I^2(x)) \quad \mbox{ as }x\to\infty.
\end{equation}
Assume $b\in(a,2a)$ and consider, for example, a service
time distribution with the Weibull integrated tail
$\overline B_I(x) = e^{-x^\beta}$, $\beta\in(0,1)$.
Then (\ref{ppp}) holds if $\bigl(\frac{b}{b-a}\bigr)^\beta>2$.

\proof\ of Theorem \ref{th.2.min.reg}.
Since $B_I\in{\cal IRV}$, both the lower bound
in Theorem \ref{th.2.min.lower} and the upper bound
in Theorem \ref{min_upper} are of the same order,
\begin{eqnarray}\label{2x.frac}
\overline B_I(2x) &=&
O\Bigl(\overline B_I\Bigl(\frac{b}{b-a}x\Bigr)\Bigr).
\end{eqnarray}

We use the notation from the previous Section. In particular, we
fix $\varepsilon>0$ and choose $L$ satisfying (\ref{LL}).
For any constant $c\ge0$, (\ref{first}) implies
\begin{eqnarray*}
\bigcup_{i=0}^{L-1} \{\sigma_{kL+i}>x+La+(L-i)c\}
&\subseteq& \bigg\{\sum_{i=0}^{L-1}
\sigma_{kL+i}-La>x\bigg\}
\subseteq \{\widehat\sigma_k>x\}.
\end{eqnarray*}
Therefore, from (\ref{first}) and (\ref{second}),
\begin{equation}\label{oo}
{\bf P}\Bigl\{\{\widehat\sigma_k>x\} \setminus
\bigcup_{i=0}^{L-1}\{\sigma_{kL+i}>x+La+(L-i)c\}\Bigr\}
= o(\overline B(x)).
\end{equation}
Take $c=(\widehat a-\widehat b)/L$. By (\ref{Z.nL.W}),
\begin{eqnarray*}
{\bf P}\{W>x\} &=&
\lim_{n\to\infty} {\bf P}\{W_{nL,1}>x\}
= \lim_{n\to\infty} {\bf P}\{W_{nL,1}>x,\widehat W_n>2x-2aL\}.
\end{eqnarray*}
Standard arguments concerning how large deviations in the single
server queue $\widehat W_n$ occur imply the relation
\begin{eqnarray*}
{\bf P}\{W>x\} &=& \lim_{n\to\infty}
\sum_{k=0}^{n-1} {\bf P}\{W_{nL,1}>x,
\widehat\sigma_k > 2x+(n-k)(\widehat a-\widehat b)\}
+o(\overline B_I(2x))\\
&=& \lim_{n\to\infty}
\sum_{i=0}^{nL-1} {\bf P}\{W_{nL,1}>x,\sigma_i>2x+(n-i)c\}
+o(\overline B_I(2x)),
\end{eqnarray*}
by (\ref{oo}). Now it follows from (\ref{LL}) that
\begin{eqnarray*}
\lefteqn{{\bf P}\{W>x\}
\le \lim_{n\to\infty} \sum_{i=0}^{nL-1}
{\bf P}\{W_{nL,1}>x,\sigma_i>2x{+}(n{-}i)(2a{-}b{-}\varepsilon)\}+
o(\overline B_I(2x))}\\
&\le& \lim_{n\to\infty} \sum_{j=1}^{nL}
{\bf P}\{W_{nL,1}>x,\sigma_{nL-j}>2x{+}j(2a{-}b{+}\varepsilon)\}
+\varepsilon O(\overline B_I(2x))+o(\overline B_I(2x))\\
&=& \lim_{n\to\infty} \Biggl(\sum_{j=1}^{N(1-\varepsilon)}
+\sum_{j=N(1-\varepsilon)}^{nL}\Biggr)+
\varepsilon O(\overline B_I(2x))
\equiv \lim_{n\to\infty}(\Sigma_1+\Sigma_2)
+\varepsilon O(\overline B_I(2x)),
\end{eqnarray*}
where $N=x/(b-a)$.
The second term admits the following estimate
\begin{eqnarray*}
\Sigma_2 &\le& \sum_{j=N(1-\varepsilon)}^\infty
{\bf P}\{\sigma>2x+j(2a-b)\}\\
&\sim& \frac{1}{2a-b} \overline B_I(2x+N(1-\varepsilon)(2a-b))
=\frac{1}{2a-b} \overline B_I\Bigl(\frac{b}{b-a}x
-\varepsilon\frac{2a-b}{b-a}x\Bigr).
\end{eqnarray*}
It follows from $B_I\in {\cal IRV}$ that,
for any $\delta >0$, there exists $\varepsilon>0$
such that
\begin{eqnarray*}
\Sigma_2 &\le&
\frac{1}{2a-b} \overline B_I\Bigl(\frac{b}{b-a}x\Bigr)
+\delta\overline B_I(2x),
\end{eqnarray*}
which coincides with the lower bound in Theorem \ref{th.2.min.lower}.

Now consider the first term $\Sigma_1$.
Since the queue is stable, one can choose $K>0$ such that
${\bf P}\{W_{n,2}\le K\}\ge 1-\varepsilon$ for all $k$. Then
\begin{eqnarray*}
\Sigma_1 &\le&
\sum_{j=1}^{N(1-\varepsilon)}
{\bf P}\{W_{nL-j,2}>K,\sigma_{nL-j}>2x+(2a-b+\varepsilon)j\}\\
&& + \sum_{j=1}^{N(1-\varepsilon)}
{\bf P}\{W_{nL,1}>x,W_{nL-j,2}\le K,
\sigma_{nL-j}>2x+(2a-b+\varepsilon)j\}\\
&\equiv& \Sigma_{1,1}+\Sigma_{1,2}.
\end{eqnarray*}
We have
\begin{eqnarray*}
\Sigma_{1,1} &=& \sum_{j=1}^{N(1-\varepsilon)}
{\bf P}\{W_{nL-j,2}>K\}{\bf P}\{\sigma_1>2x+(2a-b+\varepsilon)j\}\\
&\le& \varepsilon \sum_{j=1}^\infty
{\bf P}\{\sigma_1>2x+(2a-b)j\}
\le \frac{\varepsilon}{2a-b} \overline B_I(2x).
\end{eqnarray*}

Note that if $W_{nL-j,2}\le K$,
then $W_{nL,1}\le K+M_{nL-j+1, j-1}$. Therefore,
\begin{eqnarray*}
\Sigma_{1,2} &\le&
\sum_{j=1}^{N(1-\varepsilon)} {\bf P}\{\sigma_{nL-j}>2x+(2a-b)j,
K+M_{nL-j+1,j-1}>x\}\\
&=& \sum_{j=1}^{N(1-\varepsilon)} {\bf P}\{\sigma_{nL-j}>2x+(2a-b)j\}
{\bf P}\{K+M_{0,j-1}>x\}.
\end{eqnarray*}
Since the sequence $M_{0,j}$ stochastically increases,
\begin{eqnarray*}
\Sigma_{1,2} &\le&
{\bf P}\{K+M_{0,N(1-\varepsilon)}>x\}
\sum_{j=1}^\infty
{\bf P}\{\sigma_1>2x+(2a-b)j\}\\
&\le& {\bf P}\{M_{0,N(1-\varepsilon)}>x-K\}
\frac{1}{2a-b} \overline B_I(2x).
\end{eqnarray*}
Since
$$
\frac{x-K}{N(1-\varepsilon)}
\to \frac{b-a}{1-\varepsilon} > b-a\quad\mbox{ as }
x\to\infty,
$$
we have by (\ref{SSL})
$$
{\bf P}\{M_{0,N(1-\varepsilon)}>x-K\}
= {\bf P}\Bigl\{\frac{M_{N(1-\varepsilon)}}{N(1-\varepsilon)}
> \frac{x-K}{N(1-\varepsilon)}\Bigr\} \to 0.
$$

Thus, we have shown that the upper bound for
${\bf P}\{W>x\}$ is not bigger than the lower bound
in Theorem \ref{th.2.min.lower} plus a term of order
$$
(\varepsilon+\delta)O(B_I(2x))
\le (\varepsilon+\delta)
O\Bigl(B_I\Bigl(\frac{bx}{b-a}\Bigr)\Bigr)
$$
due to (\ref{2x.frac}). Since $\varepsilon>0$
and $\delta>0$ may be chosen as small we please,
the proof of Theorem \ref{th.2.min.reg} is complete.

\mysection{Tail asymptotics for the
two-dimensional workload vector}\label{workload}

Denote by $W^0=(W_1^0,W_2^0)$ a weak limit for the
vectors $W_n$ as $n\to\infty$. Clearly, $W=W_1^0$.

{\bf \ref{workload}.1. Maximal stability case.}
First, we obtain simple lower and upper bounds
which are equivalent up to some constant. Second, we
give (without a proof) a result related to the exact
asymptotics.

\begin{Theorem}\label{bounds.joint.max}
Let $b<a$ and $B_I\in{\mathcal L}$.
Then, as $x$, $y\to\infty$, $x\le y$,
\begin{eqnarray*}
{\bf P}\{W_1^0>x, W_2^0>y\}
&\ge& \frac{1+o(1)}{a^2}\overline B_I(x)\overline B_I(y).
\end{eqnarray*}
If, in addition, $B_I\in{\mathcal S}$, then
\begin{eqnarray*}
{\bf P}\{W_1^0>x, W_2^0>y\}
&\le& \frac{2+o(1)}{(a-b)^2}\overline B_I(x)\overline B_I(y).
\end{eqnarray*}
\end{Theorem}

\proof. Fix $\varepsilon>0$ and put $a'=a+\varepsilon$.
For $k$, $l\le n$, $k\neq l$, define the events $A_{nkl}$
and $C_{nkl}$ by the equalities
\begin{eqnarray*}
A_{nkl} &=&
\Bigl\{\sigma_k>x+(n-k)a',\ \sigma_l>y+(n-l)a'\Bigr\}
\end{eqnarray*}
and
\begin{eqnarray*}
C_{nkl} &=& \bigcap_{\stackrel{j=1}{j\ne k,l}}^n\Bigl\{
\sigma_j\le x+(n-j)a'\Bigr\}.
\end{eqnarray*}
Note that the events $A_{nkl}\cap C_{nkl}$
are disjoint for different pairs $(k,l)$ and
\begin{eqnarray*}
{\bf P}\{W_{n1}>x,W_{n2}>y\}
&\ge& \sum_{k=1}^n \sum_{l=k+1}^n
{\bf P}\{W_{n1}>x,W_{n2}>y,A_{nkl},C_{nkl}\}.
\end{eqnarray*}
Then the same calculations as in Subsection
\ref{sec.2.max.lower}.3 imply the estimate, as $x$, $y\to\infty$,
\begin{eqnarray*}
{\bf P}\{W_{n1}>x,W_{n2}>y\}
&\ge& (1+o(1))\sum_{k=1}^{n-1} \sum_{l=k+1}^{n-1}
\overline B(x+(n-k)a')\overline B(y+(n-l)a')\\
&=& (1+o(1))\sum_{k=1}^{n-1} \sum_{l=1}^{n-k-1}
\overline B(x+ka')\overline B(y+la').
\end{eqnarray*}
Hence,
\begin{eqnarray*}
{\bf P}\{W_1^0>x,W_2^0>y\}
&\ge& (1+o(1))\sum_{k=1}^\infty \sum_{l=1}^\infty
\overline B(x+ka')\overline B(y+la')
\sim \overline B_I(x)\overline B_I(y)/a'^2
\end{eqnarray*}
and the lower bound is proved.

Proceed to the upper bound. Due to construction
of the majorant $(W_n^{(1)},W_n^{(2)})$
in Section \ref{sec.2.max.upper}, we have the inequality
\begin{eqnarray*}
{\bf P}\{W_1^0>x, W_2^0>y\}
&\le& \lim_{n\to\infty}\Bigl[
{\bf P}\{W_n^{(1)}>x, W_n^{(2)}>y\}
+{\bf P}\{W_n^{(1)}>y, W_n^{(2)}>x\}\Bigr]\\
&=& 2\lim_{n\to\infty}
{\bf P}\{W_n^{(1)}>x\}{\bf P}\{W_n^{(2)}>y\}.
\end{eqnarray*}
Together with (\ref{W1.sim.BI}) it implies the desired
upper bound. Theorem \ref{bounds.joint.max} is proved.

Turn now to the exact asymptotics.
Below is the result. The proof is rather
complicated and will be presented in another paper.
Denote
$$
R(x,y) = \overline B_I(x)\overline B_I(y) + b
\int_0^\infty \overline B_I(y+za)\overline B(x+x(a-b)) dz.
$$
Recall that Theorem 1 states that
${\bf P}\{W_1^0>x\} \sim R(x,x)/a(2a-b)$
given $B_I\in {\mathcal S}$.

\begin{Theorem} \label{joint.max}
Assume $b<a$ and $B_I\in {\mathcal S}$.
Let $x, y\to \infty$, $x\le y$. Then
\begin{eqnarray*}
{\bf P}\{W_1^0>x, W_2^0>y\} &\sim&
\frac{1}{a(2a-b)}R(y,y)+\frac{1}{a^2}(R(x,y)-R(y,y)).
\end{eqnarray*}
\end{Theorem}

{\bf \ref{workload}.2. Minimal stability case.}
We prove the following

\begin{Theorem} \label{joint.min}
Assume $a<b<2a$, $B\in {\mathcal S}$, and $B_I\in {\cal IRV}$.
Let $x, y\to \infty$ in such a way that
$y/x\to c \in [1,\infty ]$. Then
\begin{eqnarray*}
{\bf P}\{W_1^0>x, W_2^0>y\}
&\sim& \frac{1}{a} \overline B_I
\Bigl(y\Bigl(1+\frac{a}{c(b-a)}\Bigr)\Bigr)
+\frac{b-a}{a(2a-b)}\overline B_I \Bigl(y\frac{b}{b-a}\Bigr).
\end{eqnarray*}
\end{Theorem}

\proof. Start with the case $c=\infty$.
From Theorem 10 in [\ref{BaF}], one can get the following:
\begin{Corollary} \label{maxx}
Assume $b\in (a,2a)$. If $B\in {\mathcal S}$
and $B_I\in {\mathcal S}$, then, as $y\to\infty$,
\begin{eqnarray*}
{\bf P}\{W_2^0>y\} &\sim&
\frac{1}{a}\overline B_I (y)
+\frac{b-a}{a(2a-b)} \overline B_I \Bigl(y\frac{b}{b-a}\Bigr).
\end{eqnarray*}
\end{Corollary}

It is clear that
$$
{\bf P}\{W_1^0>x, W_2^0>y\} \le {\bf P}\{W_2^0>y\}.
$$
On the other hand, for any $N=1$, 2, \ldots,
\begin{eqnarray*}
{\bf P}\{W_1^0>x, W_2^0>y\} &=&
\lim_{n\to\infty} {\bf P}\{W_{n,1}>x,W_{n,2}>y\}\\
&\ge& \lim_{n\to\infty} {\bf P}\Bigl\{W_{n-N,2}>y+Na,
\sum_{j=n-N}^{n-1} (\sigma_j-\tau_j)>x\Bigr\}\\
&=&
\lim_{n\to\infty} {\bf P}\{W_{n-N,2}>y+Na\}
{\bf P}\Bigl\{\sum_{j=1}^N (\sigma_j-\tau_j)>x\Bigr\}.
\end{eqnarray*}
Fix $\varepsilon>0$. Put
$N=N(x)=x(1+\varepsilon)/(b-a)$.
Then by LLN
$$
{\bf P}\Bigl\{\sum_{j=1}^N(\sigma_j-\tau_j)>x \Bigr\}
\ge 1-\varepsilon
$$
for all sufficiently large $x$ and, as $n\to\infty$,
$$
{\bf P}\{W_{n-N,2}>y+Na\} \to {\bf P}\{W_2^0>y+Na\}.
$$
Since $B_I\in {\cal IRV}$,
$$
{\bf P}\{W_2^0>y+Na\}\sim {\bf P}\{W_2^0>y\}
\quad \mbox{as} \quad  y\to\infty.
$$
By letting $\varepsilon \to 0$, we get the result.

Now consider the case $c<\infty$. If $c=1$,
then the result follows from Theorem \ref{th.2.min.reg}.
Let $c\in(1,\infty)$. We give here only a sketch of
the proof, by making links to the proof of Theorem
\ref{th.2.min.reg}.

Since
$$
{\bf P}\{W_1^0>y\} \le {\bf P}\{W_1^0>x, W_2^0>y\}
\le {\bf P}\{W_1^0>x\}
$$
and
$$
{\bf P}\{W_1^0>y\} \sim {\bf P}\{W_1^0>cx\}
\ge (K+o(1)){\bf P}\{W_1^0>x\}
$$
where $K=\inf_t \overline B_I(ct)/\overline B_I(t)>0$,
one can get from the proof of Theorem \ref{th.2.min.reg}
the following equivalences: for $N_x=x/(b-a)$,
$N_y=y/(b-a)$, and for $\varepsilon\in(0,1-1/\sqrt c)$,
\begin{eqnarray*}
\lefteqn{{\bf P}\{W_1^0>x, W_2^0>y\}}\\
&=& \lim_{n\to\infty}
\sum_{i=1}^{n-N_x(1-\varepsilon)}
{\bf P}\{W_{n,1}>x, W_{n,2}>y, \sigma_i > 2x+(n{-}i)(2a{-}b)\}
+\varepsilon O(\overline B_I(2x))\\
&=& \lim_{n\to\infty}\Biggl(
\sum_{i=1}^{n-N_y(1+\varepsilon)}
+\sum_{i=n-N_y(1+\varepsilon)}^{n-N_x(1-\varepsilon)}\Biggr)
+\varepsilon O(\overline B_I(2x))
\equiv (\Sigma_1+\Sigma_2)
+\varepsilon O(\overline B_I(2x)).
\end{eqnarray*}
Choose $K>0$ such that
${\bf P}\{W_{n,2}>K\}\le\varepsilon$ for all $n$.
Then
\begin{eqnarray*}
\Sigma_2 &=& \lim_{n\to\infty}
\sum_{i=n-N_y(1+\varepsilon)}^{n-N_x(1-\varepsilon)}
{\bf P}\{W_{i,2}\le K, W_{n,1}>x, W_{n,2}>y,
\sigma_i >2x+(n{-}i)(2a{-}b)\}\\
&&\hspace{100mm} + \varepsilon O(\overline B_I(2x)).
\end{eqnarray*}
From Lemma 2 and its Corollaries,
\begin{eqnarray*}
\Sigma_2 &=&
(1+o(1))\sum_{j=N_x(1-\varepsilon )}^{N_y(1+\varepsilon )}
{\bf P}\{\sigma_1>y+ja\}+\varepsilon O(\overline B_I(x))\\
&=& \frac{1+o(1)}{a}\Bigl(\overline B_I
\Bigl(y +\frac{x(1-\varepsilon)a}{b-a} \Bigr)
- \overline B_I \Bigl(\frac{y(b+\varepsilon a)}{b-a}\Bigr)
\Bigr)+\varepsilon O(\overline B_I(x))\\
&=& \frac{1+o(1)}{a}\Bigl(\overline B_I
\Bigl(y\Bigl(1+\frac{a}{c(b-a)}\Bigr)\Bigr)
- \overline B_I\Bigl(\frac{yb}{b-a}\Bigr)\Bigr)
+ (\varepsilon+\delta) O(\overline B_I(x)),
\end{eqnarray*}
due to $B_I\in{\cal IRV}$.
From Lemma 3 and its Corollaries, one can also conclude that, for
$i< n-N_y(1+\varepsilon )$,
if $\sigma_i<2y+(n-i)(2a-b-\varepsilon)$
and $W_{i,2} \le K$, then, with probability close to one,
both coordinates of the vector $(W_{n,1},W_{n,2})$
take values less then $y$ for all sufficiently
large $n$. From the other side, if
$\sigma_i > 2y+ (n-i) (2a-b+\varepsilon )$,
then, with probability close to one,
$y < W_{n,1} \le W_{n,2}$.
Therefore,
\begin{eqnarray*}
\Sigma_1 &=&
(1+o(1)) \sum_{j=N_y(1+\varepsilon)}^\infty
{\bf P}\{\sigma_1>2y+j(2a-b)\} + \varepsilon O(\overline B_I(x))\\
&=& \frac{1+o(1)}{2a-b}\overline B_I
\Bigl(\frac{yb}{b-a}\Bigr) + \varepsilon O(\overline B_I(x)).
\end{eqnarray*}
Summing up the terms and letting $\varepsilon$ and $\delta  \to 0$
concludes the proof.

\mysection{Comments on stationary queue length}\label{stat.q.l}

Let $Q_n$ be a queue length viewed by an arriving customer $n$, and
$Q$ its stationary version in  discrete time
(i.e. Palm-stationary). Due to the distributional Little's law,
$$
{\bf P}\{Q>n\} = {\bf P}\{W>T_n\}
$$
where $W$ is the stationary waiting time,
$T_n=\tau_1+\ldots+\tau_n$, and $W$ and $T_n$
do not depend on each other.
When a distribution of $W$ is long-tailed,
the asymptotics for ${\bf P}\{W>T_n\}$, $n\to\infty$,
have been found in [\ref{AKS}] and in [\ref{FK}].
If, in addition, $\tau_n$ has a non-lattice distribution,
there exists a stationary distribution $G$ for a queue
length in continuous time. Then, from Lemma 1 in [\ref{FK}],
$$
\overline G(n) \sim {\bf P}\{Q>n\}
\quad \mbox{as} \quad n\to\infty.
$$

\section*{\normalsize Acknowledgment}

The authors gratefully acknowledge
helpful discussions with Onno Boxma and Bert Zwart, and comments
from Daryl Daley.

\section*{\normalsize References}

\newcounter{bibcoun}
\begin{list}{[\arabic{bibcoun}]}{\usecounter{bibcoun}\itemsep=0pt}
\small

\item\label{APQ}
   S. Asmussen,
   {\it Applied Probability and Queues},
   2nd ed. (Springer, New York, 2003).

\item\label{AKS}
   S. Asmussen, C. Kluppelberg and K. Sigman,
   Sampling at subexponential times, with queueing applications,
   Stoch. Process. Appl. 79 (1999) 265--286.

\item\label{BaF} F. Baccelli and S. Foss,
   Moments and tails in monotone-separable stochastic networks,
   Ann. Appl. Probab. 14 (2004) 612--650.

\item\label{BMumn} A. A. Borovkov and A. A. Mogul'skii,
    Large deviations for Markov chains in the positive quadrant,
    Russ. Math. Surv. 56 (2001) 803--916.

\item\label{BMZ} S. Borst, M. Mandjes and A. P. Zwart,
   Exact asymptotics for fluid queues fed by heavy-tailed
   On-Off flows, Ann. Appl. Probab. 14 (2004) 903--957.

\item\label{BoF} O. J. Boxma, S. G. Foss, J.-M. Lasgouttes and
   R. Nunez Queija,
   Waiting time asymptotics in the single server queue with service
   in random order,
   Queueing Systems 46 (2004) 35--73.

\item\label{BDZ}
   O. J. Boxma, Q. Deng and A. P. Zwart,
   Waiting-time asymptotics for the $M/G/2$ queue
   with heterogeneous servers,
   Queueing Systems 40 (2002) 5--31.

\item\label{Cr}
   H. Cram\'er,
   {\it Collective risk theory}
   (Esselte, Stockholm, 1955).

\item\label{DFK}
   D. Denisov, S. Foss and D. Korshunov,
   Tail asymptotics for the supremum of a random walk
   when the mean is not finite,
   Queueing Systems 46 (2004) 15--33.

\item\label{FCh}
   S. G. Foss and N. I. Chernova,
   On optimality of FCFS discipline in multi-channel queueing systems
   and networks,
   Siberian Math. J. 42 (2001) 372--385.

\item\label{FK}
   S. Foss and D. Korshunov,
   Sampling at a random time with a heavy-tailed distribution,
   Markov Processes and Related Fields 6 (2000) 643--658.

\item\label{IMSh}
   I. A. Ignatyuk, V. Malyshev and V. Scherbakov,
   Boundary effects in large deviation problems,
   Russ. Math. Surv. 49 (1994) 41--99.

\item\label{KW55}
   J. Kiefer and J. Wolfowitz,
   On the theory of queues with many servers,
   Tran. Amer. Math. Soc. 78 (1955) 1--18.

\item\label{K}
   D.~Korshunov,
   On distribution tail of the maximum of a random walk,
   Stochastic Process. Appl. 72 (1997) 97--103.

\item\label{P}
   A. G. Pakes, On the~tails of waiting-time distribution,
   J. Appl. Probab. 12 (1975) 555--564.

\item\label{S}
   A. Scheller-Wolf,
   Further delay moment results for FIFO multiserver queues,
   Queueing Systems 34 (2000) 387--400.

\item\label{SS}
   A. Scheller-Wolf and K. Sigman,
   Delay moments for FIFO $GI/GI/s$ queues,
   Queueing Systems 25 (1997) 77--95.

\item\label{Ver}
   N. Veraverbeke,
   Asymptotic behavior of Wiener-Hopf factors
   of a random walk,
   Stochastic Process. Appl. 5 (1977) 27--37.

\item\label{W}
   W. Whitt,
   The impact of a heavy-tailed service-time
   distribution upon the $M/GI/s$ waiting-time distribution,
   Queueing Systems 36 (2000) 71--87.
\end{list}

\end{document}